\documentclass[a4paper,11pt,reqno]{amsart}

\usepackage[margin=0.9in]{geometry}
\usepackage[T1]{fontenc}
\usepackage[utf8]{inputenc}
\usepackage{lmodern}
\input{glyphtounicode}
\InputIfFileExists{glyphtounicode-cmex.tex}{}{}
\usepackage{amsmath,amssymb,amsfonts,amsthm,mathtools}
\numberwithin{equation}{section}
\usepackage{mathrsfs}
\usepackage{bm}

\usepackage[final]{graphicx}
\usepackage{epstopdf}
\usepackage{subcaption}
\usepackage{float}
\usepackage{placeins}
\usepackage{needspace}
\usepackage{algorithm}
\usepackage{algpseudocode}
\usepackage{booktabs}
\usepackage{multirow}
\usepackage{array}

\usepackage[shortlabels,inline]{enumitem}

\setlist[enumerate,1]{
  label=\textup{(\Alph*)},
  ref=\textup{(\Alph*)},
  leftmargin=*,
  itemsep=4pt,
  topsep=4pt
}

\setlist[enumerate,2]{
  label=\textup{(\roman*)},
  ref=\textup{(\roman*)},
  leftmargin=*,
  itemsep=3pt,
  topsep=3pt
}

\setlist[itemize]{
  leftmargin=*,
  itemsep=3pt,
  topsep=3pt
}

\usepackage[titletoc]{appendix}
\usepackage{setspace}

\usepackage{xcolor}

\newenvironment{revision}{\begingroup\color{black}%
  \hypersetup{linkcolor=red,citecolor=red,urlcolor=red}}{\endgroup}

\usepackage{cite}

\usepackage[
  bookmarks,
  bookmarksnumbered,
  bookmarksopen,
  breaklinks,
  linktocpage,
  colorlinks=true,
  linkcolor=blue,
  citecolor=blue,
  urlcolor=blue
]{hyperref}

\newcommand{\safeincludegraphics}[2][]{%
  \IfFileExists{#2}{%
    \includegraphics[#1]{#2}%
  }{%
    \fbox{\parbox{0.9\linewidth}{\centering Missing figure file: \texttt{\detokenize{#2}}}}%
  }%
}

\newcolumntype{L}[1]{>{\raggedright\let\newline\\\arraybackslash\hspace{0pt}}m{#1}}
\newcolumntype{C}[1]{>{\centering\let\newline\\\arraybackslash\hspace{0pt}}m{#1}}
\newcolumntype{R}[1]{>{\raggedleft\let\newline\\\arraybackslash\hspace{0pt}}m{#1}}

\theoremstyle{plain}
\newtheorem{theorem}{Theorem}[section]

\newtheorem{assumption}[theorem]{Assumption}

\theoremstyle{definition}
\newtheorem{definition}[theorem]{Definition}

\theoremstyle{remark}
\newtheorem{remark}[theorem]{Remark}

\newcommand{\ii}{\mathrm i}

\newcommand{\OmL}{\Omega_L}
\newcommand{\OmLW}{\Omega_{L,w}}
\newcommand{\Omobs}{\Omega_{\rm obs}}
\newcommand{\Bw}{\mathcal B_w}

\begin{document}
\title[Microlocal Open-Boundary Method]{A Microlocal Open-Boundary Method for Residual-Based Wave Solvers on Unbounded Domains}
	
	\author[Avy Soffer
	]{Avy Soffer
	}
	
	\address[Avy Soffer
	]{Mathematics Department \\
		 Rutgers University, New Brunswick, NJ 08903 USA}
	
	\email[Avy Soffer
	]{\href{mailto:soffer@math.rutgers.edu}{soffer@math.rutgers.edu}}

	\author[Nguyen Gia Hien]{Nguyen Gia Hien}
	
	\address[Nguyen Gia Hien]{Department of Mathematics \\
		Texas A\&M University, College Station, TX, 77843 USA}
	
	\email[Nguyen Gia Hien]{\href{mailto:giahien-nguyen@tamu.edu}{giahien-nguyen@tamu.edu}}
	
	\author[Minh-Binh Tran]{Minh-Binh Tran}
	\address[Minh-Binh Tran]{Department of Mathematics \\
		Texas A\&M University, College Station, TX, 77843 USA}
	
	\email[Minh-Binh Tran]{\href{mailto:minhbinh@tamu.edu}{minhbinh@tamu.edu}}
	\thanks{N. G. Hien and M.-B. T are  funded in part by     NSF CAREER  DMS-2303146, and NSF Grants DMS-2204795, DMS-2305523,  DMS-2306379. A.S. is funded in part by NSF Grant DMS-2205931.   }

\begin{abstract}
\emergencystretch=1em
\begin{revision}
We couple the established time-dependent phase-space filter (TDPSF)
with a slabwise physics-informed neural network (PINN) for open-domain
wave propagation. The neural solver approximates the evolution on each
time slab, and localized Fourier multipliers attenuate outgoing buffer
components according to their group velocity. The resulting Microlocal
PINN Filter (TDPSF-PINN) combines neural residual minimization, sequential buffer
filtering, and transfer of the computed endpoint to the next slab.

We present $3$ numerical tests: an anisotropic focusing nonlinear Schr\"odinger packet
with a PML comparator, acoustic radiation with a 
CRBC comparator, and acoustic scattering by three rigid disks with a
second-order Bayliss--Turkel comparator.   At the reported budgets, TDPSF--PINN shows better accuracy and stability  in
all 3 tests in the long-time simulations.
\end{revision}
\end{abstract}
    \maketitle

\section{Introduction}
\begin{revision}
Open-domain wave propagation is a fundamental challenge in computational
science and engineering because the physical domain is unbounded while
the numerical computation must be performed on a finite region.
Schr\"odinger equations, nonlinear Schr\"odinger equations, acoustic
waves, electromagnetic waves, elastic waves, and anisotropic hyperbolic
systems are naturally posed on \(\mathbb R^d\) or on exterior domains.
In computations, however, one must replace the unbounded domain by a
finite computational region.  This truncation introduces an artificial
boundary, and the boundary treatment can strongly affect the solution.
Outgoing waves may be reflected back into the physical region,
periodically wrapped around the computational box, or over-damped by an
absorbing layer.  These artifacts are especially damaging for dispersive
and anisotropic wave equations, where the outgoing character of a wave
is not a purely local property of its trace on the boundary.

Classical open-boundary techniques include absorbing boundary
conditions, nonreflecting boundary conditions, perfectly matched layers
(PMLs), damping or sponge layers, and rational approximations of exact
radiation operators.  For homogeneous isotropic wave equations, these
methods are  highly effective.  Absorbing and radiation boundary
conditions go back to the work of Engquist--Majda, Bayliss--Turkel, and
related developments
\cite{MR0471386,MR0436612,MR517938,MR596431,MR1819643,MR2032866,givoli1991non}.
The PML, introduced by B\'erenger for Maxwell's equations, is among the
most widely used artificial-boundary techniques
\cite{Berenger1994,MR1412240}.  It may be interpreted as a complex
coordinate stretching or, more generally, as a transformation of the
dispersion relation.

However, these methods are not uniformly reliable for anisotropic or
multi-branch wave systems.  In particular, for certain anisotropic
problems the PML can be exponentially unstable at the PDE level.  This
phenomenon has been observed for the linearized Euler equations about a mean
flow and related systems
\cite{hu:unstablePML,BecacheFauqueuxJoly2003,abarbanel:PMLinstability,
	abarbanel:PMLinstability2}.  The underlying difficulty is geometric and can be explained as follows: a
local absorbing layer may shift parts of the spectrum into an unstable
region while a scalar local absorber does not necessarily distinguish
different wave families and their group velocities.  Similarly,
high-order or optimized radiation conditions can be very effective when
the radiation operator has a tractable scalar or isotropic structure;
however their direct extension to general anisotropic or multi-branch
systems is limited.  Damping layers remain useful in practice, but they
are not selective  due to the fact that  they attenuate incoming, outgoing, glancing,
acoustic, vortical, or dispersive content according to location rather
than according to the physical propagation direction of each wave
packet.

The time-dependent phase-space filter (TDPSF) provides a different
viewpoint on open boundaries
\cite{HagstromNguyenSofferStucchioTran2026,SofferStucchio2007,
	us:multiscale,SofferStucchioTran2023}.  Instead of imposing a local boundary
condition at the edge of the computational box, TDPSF identifies
portions of the numerical solution that are localized near the boundary
and whose group velocities point outward, and removes only those
components before they can interact with the artificial boundary.  

Physics-informed neural networks (PINNs) provide a flexible
residual-minimization framework for approximating solutions of partial
differential equations from equation residuals, initial data, boundary
data, and possibly observations \cite{RaissiPerdikarisKarniadakis2019}.
Boundary and initial conditions in PINNs are commonly imposed either
weakly through penalty terms or strongly through constrained trial
spaces; see, for example,
\cite{BerroneCanutoPintoreSukumar2023,SukumarSrivastava2022,
	AlkhadhrAlmekkawy2023}.  For wave problems on truncated domains,
absorbing-boundary and PML-type treatments have also been incorporated
into PINN losses or residual formulations
\cite{RenRaoChenWangSunLiu2024,WuAghamiryOpertoMa2022,
	AbediPardoAlkhalifah2025}. These
approaches remain local or layer-based: they act through boundary
traces, boundary residuals, or damped/PML residuals, rather than through
a microlocal classification of outgoing wave packets. For
wave propagation, however, this local boundary-loss strategy inherits a
basic limitation: radiation is microlocal, depending on both position
and frequency, whereas a local boundary penalty acts only on traces or
local residuals.  A boundary loss does not directly distinguish
incoming, outgoing, glancing, and branch-dependent modes.

This distinction is especially important for anisotropic or
multi-branch systems.  In such problems, phase velocity and group
velocity need not be aligned, and different dispersion branches may
propagate in different directions at the same frequency.  These
geometric effects are also responsible for known failures of local
absorbing mechanisms As an example, PMLs can be exponentially unstable
for some anisotropic systems and for the linearized Euler equations
about a mean flow
\cite{BecacheFauqueuxJoly2003,hu:unstablePML,abarbanel:PMLinstability,
	abarbanel:PMLinstability2}.  As thus, a boundary treatment based only on
local traces, or on a naive scalar sign condition involving \(k\cdot n\),
can remove or retain the wrong components.

The purpose of this work is to incorporate the microlocal time-dependent
phase-space-filter  mechanism \cite{HagstromNguyenSofferStucchioTran2026,SofferStucchio2007,
	us:multiscale,SofferStucchioTran2023} into a PINN framework.  Let

\[
\OmL=[-L,L]^d
\]
denote the physical computational box and let

\[
\OmLW=[-L-w,L+w]^d,
\qquad
\Bw=\OmLW\setminus\OmL,
\]
denote the extended box and its buffer layer.  On each time slab

\[
I_m=[t_m,t_{m+1}],
\]
a neural network approximates the interior evolution on
\(I_m\times\OmLW\).  At the end of the slab, the network output is
evaluated on a grid, and a phase-space filter removes the outgoing
components located in the buffer.  The filtered state is then used as
the initial condition for the next slab:

\[
u_m^{\rm filt}
\xrightarrow{\text{PINN on } I_m}
u_{m+1}^{\rm raw}
\xrightarrow{\text{phase-space filter}}
u_{m+1}^{\rm filt}.
\]

The phase-space filtering component of the present work follows the
TDPSF philosophy \cite{HagstromNguyenSofferStucchioTran2026,SofferStucchio2007,
	us:multiscale,SofferStucchioTran2023}, in which outgoing radiation is not characterized by a local
boundary trace alone.  Instead, the solution is decomposed into
localized wave packets, coherent states, or localized Fourier
components and a packet near the artificial boundary is removed only if
its group velocity points outward.

The main contribution of this paper is not a new neural-network
architecture.  Rather, the contribution is the coupling of a
residual-based neural interior solver with a microlocal open-boundary
mechanism.   When a forward problem has
known constant coefficients and a spectral formulation is available, a
classical FFT, split-step, or spectral solver combined with a
time-dependent phase-space filter remains a natural and often superior
choice.  The proposed Microlocal PINN Filter (TDPSF-PINN) framework is intended for settings in
which the flexibility of a neural residual formulation is useful:
nonlinear equations, sparse or off-grid observations, data
assimilation, variable or unknown interior coefficients, and media or
geometries for which the interior evolution is not conveniently
diagonalized by FFT.  In the TDPSF-PINN framework, the phase-space filter
supplies the radiation mechanism, while the PINN supplies the flexible
interior approximation.

The paper is organized as follows.  Section~\ref{sec:geometry}
introduces the computational geometry, coherent-state notation, and
localized FFT-window filters.  Section~\ref{sec:linear-sch} develops
the linear Schr\"odinger algorithm with potential. Based on the preparation of Section~\ref{sec:linear-sch}, Section~\ref{sec:nls}
extends the linear framework to nonlinear Schr\"odinger equations and discusses
coherent-state deletion, FFT-window filtering, and nonlinear buffer
diagnostics.  Section~\ref{sec:anisotropic} treats anisotropic
multi-branch wave systems and defines the branchwise group-velocity
filter.   
The numerical study in Section~\ref{sec:results} consists of 3  problems with prescribed
coefficients and forcing. Test 1 compares TDPSF--PINN with the specified
complex-coordinate PML--PINN \cite{Zheng2007} for a nonlinear Schr\"odinger equation. Test 2 compares TDPSF--PINN with Hagstrom--Warburton's CRBC method \cite{HagstromWarburton2009} (coupled with PINN)
for zero-flow acoustic radiation. Test 3 compares TDPSF--PINN with a
second-order Bayliss--Turkel method \cite{MR596431} (also coupled with PINN in the same settings) for a multiple scattering problem on 3 rigid
disks. In those tests, each comparison uses paired seeds and shared field-network
architecture, collocation rules, and optimizer budgets within the test.
In our simulations, we compare the   behavior of the computed solutions    following the long-time
simulation approach of Hagstrom et al.~\cite{HagstromNguyenSofferStucchioTran2026}.
   At the reported
budgets, the  TDPSF--PINN method is more stable and achieves smaller  errors in all 3
tests, as time evolves.

\end{revision}

\section{Computational geometry and phase-space notation}
\label{sec:geometry}

Let \(d\geq 1\).  We use the nested boxes
$$
  \OmL=[-L,L]^d,
  \qquad
  \OmLW=[-L-w,L+w]^d,
  \qquad
  \Bw=\OmLW\setminus\OmL,
$$
where \(w>0\) is the buffer width.  The solution of physical interest
is the restriction to \(\OmL\), or to a smaller observation region
\[
  \Omobs\Subset \OmL.
\]
The enlarged box \(\OmLW\) is used only for the implementation of the
phase-space filter.

For each coordinate direction \(j\in\{1,\ldots,d\}\) and each sign
\(s\in\{+,-\}\), we define the side buffer
$$
  \mathcal B_{w,j}^{s}
  =
  \begin{cases}
  \{x\in\OmLW: x_j\in[L,L+w]\}, & s=+,\\[1mm]
  \{x\in\OmLW: x_j\in[-L-w,-L]\}, & s=-.
  \end{cases}
$$
The corresponding outward unit normal is then computed as

\[
  n_j^{s}
  =
  \begin{cases}
  e_j, & s=+,\\
  -e_j, & s=-.
  \end{cases}
\]

\subsection{Coherent-state and FFT-window phase-space filters}
\label{subsec:phase-space-filter-notation}

We first recall the transforms used in the phase-space filter.  Let
\(\mathcal F, \mathcal F^{-1}\) denote the spatial Fourier transform and its inverse
\[
  (\mathcal F u)(k)
  =
  \widehat u(k)
  =
  \int_{\mathbb R^d} e^{-i x\cdot k}u(x)\,\mathrm dx,\ \ \ 
  (\mathcal F^{-1}\widehat u)(x)
  =
  \frac{1}{(2\pi)^d}
  \int_{\mathbb R^d} e^{i x\cdot k}\widehat u(k)\,\mathrm dk.
\]
On a finite computational box these operators are implemented by FFT
and inverse FFT.

Let \(g_\sigma\) be the normalized Gaussian window
\[
  g_\sigma(x)
  =
  (\pi\sigma^2)^{-d/4}
  \exp\!\left(-\frac{|x|^2}{2\sigma^2}\right).
\]
The associated windowed Fourier transform 
is defined as follows
$$
  (\mathcal G_\sigma u)(q,p)
  =
  \int_{\mathbb R^d}
  u(x)\overline{g_\sigma(x-q)}e^{-ip\cdot x}\,\mathrm dx,
$$
in which \(q\in\mathbb R^d\) is the packet center in physical space and
\(p\in\mathbb R^d\) is the packet frequency, or momentum.

The adjoint
operator \(\mathcal G_\sigma^*\) reconstructs a function from
phase-space coefficients as follows

\[
  (\mathcal G_\sigma^* a)(x)
  =
  \int_{\mathbb R^d}\int_{\mathbb R^d}
  a(q,p)g_\sigma(x-q)e^{ip\cdot x}
  \,\mathrm dq\,\mathrm dp,
\]
up to the normalization constant determined by the chosen Fourier
convention.

In computations, \((q,p)\) is sampled on a lattice
\(
  \Lambda_q\times\Lambda_p.
\)
The reconstruction is then replaced by a discrete frame synthesis
formula, which reads

\[
  \mathcal G_\sigma^* a
  \approx
  \sum_{q\in\Lambda_q}
  \sum_{p\in\Lambda_p}
  w_{q,p}\,
  a(q,p)\,
  g_\sigma(\,\cdot-q\,)e^{ip\cdot(\cdot)},
\]
where \(w_{q,p}\) are quadrature or frame weights.   

\subsection{Outgoing masks and branchwise group velocity}
\label{subsec:outgoing-masks}

Let \(\omega_\ell(k)\) denote the dispersion relation of branch
\(\ell\).

The associated group velocity is
$$
  v_{g,\ell}(k)=\nabla_k\omega_\ell(k).
$$
For scalar equations in which there is only one branch   we write simply
$$
  \omega(k),
  \qquad
  v_g(k)=\nabla_k\omega(k).
$$

The basic principle of the filter is that a wave packet should be
removed only if it is both close to the artificial boundary and moving
out of the computational domain.  In the scalar case, this condition is
encoded by an outgoing phase-space mask defined as follows
\[
  \chi_{\rm out}(q,k)
  =
  \eta_{\rm bdry}(q)\,
  S_\delta\!\left(v_g(k)\cdot n(q)\right),
\]
in which $$\eta_{\rm bdry}$$ is a smooth spatial cutoff supported in the
boundary buffer, $$n(q)$$ is the outward normal at the boundary point
nearest to \(q\), and \(S_\delta\) is a smooth approximation of the
Heaviside function.  A typical choice is as follows

\[
  S_\delta(r)=\frac{1}{1+e^{-r/\delta}}.
\]
As thus, we find $$\chi_{\rm out}(q,k)\approx 1$$ for wave packets near the
artificial boundary whose group velocity points outward, and
$$\chi_{\rm out}(q,k)\approx 0$$ for packets that are either away from
the boundary or not outgoing.

For systems, the filtering is performed branch by branch.  Let \(\Pi_\ell(k)\) be the spectral projector associated with the branch
\(\omega_\ell(k)\).  For a simple branch, we have

\[
  \Pi_\ell(k)=d_\ell(k)d_\ell(k)^*,
\]
where \(d_\ell(k)\) is the normalized eigenvector satisfying
$$M(k)d_\ell(k)=\omega_\ell(k)d_\ell(k).$$  Near repeated eigenvalues,
\(\Pi_\ell(k)\) is replaced by the projector onto the corresponding
spectral cluster. 

\begin{remark} 
The notation \(M(k)\) is used only when the far-field dispersion relation
is obtained from a matrix-valued symbol.  

For the standard linear Schr\"odinger equation in the far field,
$$\ii\partial_t\psi=-\frac12\Delta\psi+V_\infty\psi,$$ the plane wave
ansatz \(\psi(t,x)=a e^{\ii(k\cdot x-\omega t)}\) gives us directly
$$\omega(k)=\frac12|k|^2+V_\infty.$$
Thus the scalar analogue of \(M(k)\) is simply
$M_{\rm Sch}(k)=\omega(k)=\frac12|k|^2+V_\infty$, and the spectral
projector is trivial: \(\Pi_1(k)=1\).
The outgoing condition at side \((j,s)\) is therefore
\(
  \nabla_k\omega(k)\cdot n_j^s
  =
  k\cdot n_j^s>0.
\)

For the anisotropic Schr\"odinger equation
\[
  \ii\partial_t\psi
  =
  -\frac12\nabla\cdot(A\nabla\psi)+V_\infty\psi,
  \qquad A=A^T>0,
\]
the corresponding scalar symbol is then

\[
  M_{\rm Sch,A}(k)=\omega(k)
  =
  \frac12 k^T A k+V_\infty.
\]
Again there is only one branch and \(\Pi(k)=1\).  The group velocity is
$$v_g(k)=\nabla_k\omega(k)=Ak$$  and hence the outgoing condition is
\(Ak\cdot n_j^s>0\).

For a first-order anisotropic wave system
\[
  \partial_t U
  =
  \mathcal H U,
  \qquad
  \mathcal H=\sum_{r=1}^d A_r\partial_{x_r}+B,
\]
with \(A_r=A_r^*\) and \(B^*=-B\), the Fourier symbol is
\(
  \widehat{\mathcal H}(k)
  =
  \ii\sum_{r=1}^d k_r A_r+B.
\)
Using the plane wave convention
\(
  U(t,x)=d(k)e^{\ii(k\cdot x-\omega t)},
\)
the dispersion matrix is the Hermitian matrix
\[
  M(k)
  =
  \ii\,\widehat{\mathcal H}(k)
  =
  -\sum_{r=1}^d k_r A_r+\ii B.
\]
The dispersion branches are the eigenvalues of \(M(k)\):
\[
  M(k)d_\ell(k)=\omega_\ell(k)d_\ell(k),
  \qquad \ell=1,\ldots,q.
\]
When \(\omega_\ell(k)\) is simple, the associated spectral projector is
\(
  \Pi_\ell(k)=d_\ell(k)d_\ell(k)^*.
\)
The branchwise group velocity is \(v_{g,\ell}(k)=\nabla_k\omega_\ell(k)\),
and the outgoing condition at side \((j,s)\) is
\(v_{g,\ell}(k)\cdot n_j^s>0\).
 
\end{remark}

For side \((j,s)\), branch \(\ell\) is outgoing when
\(v_{g,\ell}(k)\cdot n_j^s>0\), and we denote the corresponding outgoing
Fourier multiplier by \(P_{\ell,j}^{s}(k)\).  In the scalar case, we have
$$P_j^s(k)=S_\delta(v_g(k)\cdot n_j^s),$$ while for a system a typical
branchwise multiplier is
\(P_{\ell,j}^{s}(k)=S_\delta(v_{g,\ell}(k)\cdot n_j^s)\Pi_\ell(k)\).
Thus \(P_{\ell,j}^{s}(k)\) selects the branch-\(\ell\) component whose
group velocity points outward through side \((j,s)\).

\subsection{Coherent-state filter}
\label{subsec:coherent-state-filter}

With the notations above, the coherent-state version of the outgoing
filter reads as follows
\begin{equation}\label{eq:gabor-filter}
  \mathsf F u
  =
  \mathcal G_\sigma^*
  \big[(1-\chi_{\rm out})\mathcal G_\sigma u\big].
\end{equation}
Equivalently, one first decomposes \(u\) into localized wave packets,
then multiplies the phase-space coefficients by \(1-\chi_{\rm out}\),
and finally reconstructs the solution.  Wave packets classified as
outgoing are deleted, while interior, incoming, and tangential packets
are retained in the computations.

\subsection{Localized FFT-window filter}
\label{subsec:fft-window-filter}

For higher-dimensional computations and for anisotropic systems, we use
a localized FFT-window version of the same principle.  Let
\(\eta_j^s\) be a smooth spatial cutoff supported in the side buffer
\(\mathcal B_{w,j}^{s}\):
$$
  0\leq \eta_j^s\leq 1,
  \qquad
  \operatorname{supp}\eta_j^s\subset \mathcal B_{w,j}^{s}.
$$
If branch-dependent buffers or branch-dependent cutoffs are needed, we
write \(\eta_{\ell,j}^{s}\) and
\(\mathcal B_{w,\ell,j}^{s}\).  Otherwise the branch index is suppressed
on the spatial cutoff.

For a scalar equation, the localized outgoing component through side
\((j,s)\) is
\[
  \eta_j^s(x)\mathcal F^{-1}
  \left[
    P_j^s(k)\mathcal F(\eta_j^s u)(k)
  \right].
\]
The corresponding FFT-window filter is
\begin{equation}\label{eq:simplified-filter-template-scalar}
  u
  \longmapsto
  u-
  \eta_j^s(x)\mathcal F^{-1}
  \left[
    P_j^s(k)\mathcal F(\eta_j^s u)(k)
  \right].
\end{equation}
In the case of systems, the branchwise version is
\begin{equation}\label{eq:simplified-filter-template-system}
  u
  \longmapsto
  u-
  \eta_{\ell,j}^{s}(x)\mathcal F^{-1}
  \left[
    P_{\ell,j}^{s}(k)\mathcal F(\eta_{\ell,j}^{s}u)(k)
  \right].
\end{equation}
This should be interpreted as follows.     First, multiplication by
\(\eta_{\ell,j}^{s}\) localizes the solution to the buffer layer near
the side    \((j,s)\).  Second, the Fourier transform decomposes this
localized piece into frequencies.    Third, \(P_{\ell,j}^{s}(k)\) selects
the outgoing branch-\(\ell\) components whose group velocity points
through that side.  Fourth, the inverse Fourier transform reconstructs
the outgoing component in physical space.  Finally, this component is
multiplied again by \(\eta_{\ell,j}^{s}\) and subtracted from \(u\).

Applying this operation over all sides and branches gives the schematic
filter
\begin{equation}\label{eq:fft-filter-all-sides-branches}
  \mathsf F u
  =
  u-
  \sum_{j=1}^d
  \sum_{s\in\{+,-\}}
  \sum_{\ell}
  \eta_{\ell,j}^{s}(x)\mathcal F^{-1}
  \left[
    P_{\ell,j}^{s}(k)\mathcal F(\eta_{\ell,j}^{s}u)(k)
  \right].
\end{equation}
If the spatial cutoffs are independent of the branch, one simply sets
$$\eta_{\ell,j}^{s}=\eta_j^s.$$

\subsection{Construction of the side windows}
\label{subsec:side-window-construction}

We now describe a concrete construction of the cutoff \(\eta_j^s\).    Let
\[
  G_\sigma(x)
  =
  \frac{1}{(\sigma\sqrt{\pi})^d}
  \exp\!\left(-\frac{|x|^2}{\sigma^2}\right),
\]
and let \(I_j^s\) be the indicator of a strip inside
\(\mathcal B_{w,j}^{s}\).  For example, for the positive side normal to
the \(x_j\)-direction, set
\[
  I_j^+(x)
  =
  \mathbf 1_{\mathcal S_j^+}(x),
\]
where

\[
  \mathcal S_j^+
  =
  \left\{
  x:
  x_j\in [L+w/3,L+2w/3],
  \quad
  x_m\in[-L-2w/3,L+2w/3]\ \text{for }m\neq j
  \right\}.
\]
For the negative side, set

\[
  I_j^-(x)
  =
  \mathbf 1_{\mathcal S_j^-}(x),
\]
where
\[
  \mathcal S_j^-
  =
  \left\{
  x:
  x_j\in [-L-2w/3,-L-w/3],
  \quad
  x_m\in[-L-2w/3,L+2w/3]\ \text{for }m\neq j
  \right\}.
\]
We then define
\[
  \eta_j^s=G_\sigma*I_j^s.
\]
If necessary, the resulting function is smoothly truncated near corners
or renormalized so that \(0\leq \eta_j^s\leq 1\).  For systems, one may
use the same window for all branches, $$\eta_{\ell,j}^{s}=\eta_j^s,$$ or
introduce branch-dependent windows \(\eta_{\ell,j}^{s}\) if different
branches require different buffer widths.

T In practice, one
chooses \(\sigma\), a frequency buffer \(k_b\), and a tolerance
\(\delta\) so that
\begin{equation}\label{eq:sigma-window-bounds}
  C_1 k_b^{-1}\sqrt{\log(\delta^{-1})}
  \leq
  \sigma
  \leq
  C_2 \frac{w}{\sqrt{\log(\delta^{-1})}}.
\end{equation}
  Consequently one obtains the expected uncertainty scaling
\begin{equation}\label{eq:w-scaling}
  w\gtrsim k_b^{-1}\log(\delta^{-1}).
\end{equation}
The constants in \eqref{eq:sigma-window-bounds} depend on the precise
definition of the window, the spatial dimension, and the Fourier
transform convention.

\section{Linear Schr\"odinger equation with potential}
\label{sec:linear-sch}

We first consider the linear Schr\"odinger equation
\begin{equation}\label{eq:schrodinger-full}
  \ii\partial_t\psi(t,x)
  =
  -\frac12\Delta\psi(t,x)+V(x)\psi(t,x)+f(t,x),
  \qquad x\in\mathbb R^d,
\end{equation}
with initial condition
\begin{equation}\label{eq:schrodinger-ic}
  \psi(0,x)=\psi_0(x).
\end{equation}
The potential \(V\) is assumed to be real-valued.  The forcing term
\(f\) may be set to zero.

\begin{assumption} 
\label{ass:far-field-linear}
There exists a constant \(V_\infty\in\mathbb R\) such that
\[
  V(x)=V_\infty
  \qquad \text{for }x\in\Bw.
\]
More generally, the algorithm remains meaningful if \(V-V_\infty\) and
\(\nabla V\) are small in the buffer region \(\Bw\).  If the potential
is long-range or strongly varying in the buffer, then the outgoing
classifier should be replaced by a local Hamiltonian or ray-tracing
classifier.
\end{assumption}

Under Assumption~\ref{ass:far-field-linear}, the far-field dispersion
relation is
\begin{equation}\label{eq:schrodinger-dispersion}
  \omega(k)=\frac12|k|^2+V_\infty,
  \qquad
  v_g(k)=\nabla_k\omega(k)=k.
\end{equation}
In some applications the far-field Schr\"odinger operator may have
direction-dependent dispersion.  The corresponding anisotropic far-field equation is
\[
  \ii\partial_t\psi
  =
  -\frac12\nabla\cdot(A\nabla\psi)+V_\infty\psi .
\]
For the plane wave \(\psi(t,x)=a e^{\ii(k\cdot x-\omega t)}\), we obtain
\(\omega(k)=\frac12 k^T A k+V_\infty\).
Consequently the group velocity is
\begin{equation}\label{eq:anisotropic-sch-dispersion}
  v_g(k)
  =
  \nabla_k\omega(k)
  =
  Ak.
\end{equation}
Hence, in the anisotropic scalar Schr\"odinger case, outgoing packets
must be classified using \(v_g(k)\cdot n_j^s=Ak\cdot n_j^s\), not
merely using \(k\cdot n_j^s\).  This distinction is important
because the direction of energy or wave-packet propagation need not be
parallel to the wave vector \(k\).

\subsection{Neural representation and PDE residual}
\label{subsec:linear-sch-pinn-residual}

We write the neural approximation in real and imaginary parts:
\[
  \psi_\theta(t,x)=a_\theta(t,x)+\ii b_\theta(t,x),
\]
where $$a_\theta,b_\theta:\mathbb R\times\mathbb R^d\to\mathbb R$$ are
the two outputs of the neural network.  The complex residual associated
with \eqref{eq:schrodinger-full} is
\begin{equation}\label{eq:schrodinger-complex-residual}
  \mathcal R_\theta
  =
  \ii\partial_t\psi_\theta
  +
  \frac12\Delta\psi_\theta
  -
  V\psi_\theta
  -
  f.
\end{equation}
If \(f=f_R+\ii f_I\), then the real and imaginary residuals are
\begin{align}
  \mathcal R_\theta^R
  &=
  -\partial_t b_\theta
  +\frac12\Delta a_\theta
  -Va_\theta
  -f_R,
  \label{eq:sch-residual-real}\\
  \mathcal R_\theta^I
  &=
  \partial_t a_\theta
  +\frac12\Delta b_\theta
  -Vb_\theta
  -f_I.
  \label{eq:sch-residual-imag}
\end{align}
On a time slab
\(
  I_m=[t_m,t_{m+1}],
\)
the residual loss is approximated by collocation:

\begin{equation}\label{eq:sch-residual-loss}
  \mathcal L_{\rm PDE}^{(m)}(\theta)
  =
  \frac{1}{N_r}
  \sum_{r=1}^{N_r}
  \left(
  |\mathcal R_\theta^R(\tau_r,x_r)|^2
  +
  |\mathcal R_\theta^I(\tau_r,x_r)|^2
  \right),
\end{equation}
where
\(
  (\tau_r,x_r)\in I_m\times\OmLW.
\)
The collocation points should be denser in the physical box \(\Omega_L\),
near the support of the potential, and in the side buffers
\(\mathcal B_{w,j}^{s}\), \(j=1,\ldots,d\), \(s\in\{+,-\}\).

\subsection{Outgoing frequency masks for the linear Schr\"odinger equation}
\label{subsec:linear-sch-outgoing-masks}

For the side \((j,s)\), where \(j\in\{1,\ldots,d\}\) and
\(s\in\{+,-\}\), the outward unit normal is \(n_j^s\).  The outgoing
condition for a scalar Schr\"odinger packet is

\begin{equation}\label{eq:sch-outgoing-condition}
  v_g(k)\cdot n_j^s>0.
\end{equation}
For the standard isotropic Schr\"odinger equation, \(v_g(k)=k\), so
\eqref{eq:sch-outgoing-condition} becomes \(k\cdot n_j^s>0\), or
equivalently \(s k_j>0\), where \(s=+1\) on the positive \(x_j\)-face
and \(s=-1\) on the negative \(x_j\)-face.

To avoid aggressively filtering packets close to the glancing set
\(
  v_g(k)\cdot n_j^s=0,
\)
we introduce a glancing buffer \(\gamma>0\) and define the smooth
outgoing Fourier multiplier

\begin{equation}\label{eq:sch-frequency-mask}
  P_j^s(k)
  =
  S_\alpha\!\left(v_g(k)\cdot n_j^s-\gamma\right),
\end{equation}
where
\(
  S_\alpha(r)=\frac{1}{1+\exp(-r/\alpha)}
\)
is a smoothed Heaviside function.  The parameter \(\alpha>0\) controls
the sharpness of the transition, while \(\gamma>0\) prevents packets
near the glancing region from being removed too strongly.

For the anisotropic Schr\"odinger operator
\(-\frac12\nabla\cdot(A\nabla)\), the multiplier is still given by
\eqref{eq:sch-frequency-mask}, but with the choice $$v_g(k)=Ak.$$   Thus the
correct outgoing condition is $$Ak\cdot n_j^s>0$$ not simply
\(k\cdot n_j^s>0\).

\subsection{Linear Schr\"odinger outgoing filter}
\label{subsec:linear-sch-outgoing-filter}

Let \(u:\OmLW\to\mathbb C\) be a grid function.  For the side
\((j,s)\), let \(\eta_j^s\) be the side window supported in
\(\mathcal B_{w,j}^{s}\), as constructed in
Section~\ref{subsec:side-window-construction}.  The outgoing component
through side \((j,s)\) is defined by
\begin{equation}\label{eq:sch-outgoing-component}
  \mathcal O_j^s u
  =
  \eta_j^s(x)
  \mathcal F^{-1}
  \left[
    P_j^s(k)\,
    \mathcal F(\eta_j^s u)(k)
  \right](x).
\end{equation}
The corresponding one-side filter is
\begin{equation}\label{eq:sch-filter-one-side}
  \mathsf F_j^s u
  =
  u-\mathcal O_j^s u.
\end{equation}
Applying the filter over all sides gives the complete linear
Schr\"odinger filter
\begin{equation}\label{eq:sch-complete-filter}
  \mathsf F_{\rm Sch}u
  =
  \left(
  \prod_{j=1}^{d}
  \mathsf F_j^{-}\mathsf F_j^{+}
  \right)u.
\end{equation}

\begin{algorithm}[tbp]
\caption{Precomputation for the linear Schr\"odinger phase-space filter}
\label{alg:sch-precompute}
\begin{algorithmic}[1]
\Require Dimension \(d\), physical half-width \(L\), buffer width \(w\),
tolerance \(\delta\), glancing buffer \(\gamma\), smoothing width
\(\alpha\), grid size \(N_x\), far-field potential \(V_\infty\), maximum
resolved frequency \(K_{\max}\).
\Ensure Side windows \(\eta_j^s\), outgoing masks \(P_j^s\), filter time
step \(\Delta T_{\rm filt}\).
\State Define the extended periodic box
\[
  \OmLW=[-L-w,L+w]^d.
\]
\State Build a uniform tensor grid on \(\OmLW\) and the corresponding
Fourier grid \(\{k_r\}\).
\State Construct the Gaussian-mollified side windows $$\eta_j^s$$ for
all \(j=1,\ldots,d\) and \(s\in\{+,-\}\).
\State For each Fourier grid point \(k_r\), compute the far-field group
velocity:
\[
  v_g(k_r)=k_r
\]
for the standard Schr\"odinger operator, or
\[
  v_g(k_r)=Ak_r
\]
for an anisotropic mass tensor \(A\).
\State Define the outgoing mask values
\[
   P_j^s(k_r)
   =
   S_\alpha\!\left(v_g(k_r)\cdot n_j^s-\gamma\right).
\]
\State Set
\[
   v_{\max}
   =
   \sup_{|k_r|\leq K_{\max}}|v_g(k_r)|.
\]
For the standard Schr\"odinger equation, set $$v_{\max}=K_{\max}.$$
\State Choose the filtering time step so that
\[
   \Delta T_{\rm filt}\leq \frac{w}{3v_{\max}}.
\]
\State Check the uncertainty condition
\[
  C_1 k_b^{-1}\sqrt{\log(\delta^{-1})}
  \leq
  \sigma
  \leq
  C_2\frac{w}{\sqrt{\log(\delta^{-1})}}.
\]
If this condition fails, increase \(w\), increase the frequency buffer,
or reduce the requested tolerance.
\end{algorithmic}
\end{algorithm}

\begin{algorithm}[tbp]
\caption{One application of the linear Schr\"odinger outgoing filter}
\label{alg:sch-filter}
\begin{algorithmic}[1]
\Require Grid values \(u(x_i)\) on \(\OmLW\), side windows \(\eta_j^s\),
and outgoing masks \(P_j^s\).
\Ensure Filtered grid values \(u^{\rm filt}\).
\State Set \(u^{(0)}=u\) and \(r=0\).
\For{\(j=1,\ldots,d\)}
  \For{\(s\in\{+,-\}\)}
    \State Localize the current state:
    \[
      z(x_i)=\eta_j^s(x_i)u^{(r)}(x_i).
    \]
    \State Compute $$\widehat z(k)=\mathcal F z$$ by FFT.
    \State Extract the outgoing component in frequency:
    \[
      \widehat z_{\rm out}(k)=P_j^s(k)\widehat z(k).
    \]
    \State Transform back and localize again:
    \[
      o(x_i)=\eta_j^s(x_i)\mathcal F^{-1}\widehat z_{\rm out}(x_i).
    \]
    \State Remove the outgoing component:
    \[
      u^{(r+1)}(x_i)=u^{(r)}(x_i)-o(x_i).
    \]
    \State Set$$r\gets r+1.$$
  \EndFor
\EndFor
\State Return $$u^{\rm filt}=u^{(r)}.$$
\end{algorithmic}
\end{algorithm}

\subsection{Slabwise PINN--TDPSF algorithm for the linear Schr\"odinger equation}
\label{subsec:linear-sch-slabwise-pinn}

Let
\(
  0=t_0<t_1<\cdots<t_M=T,
  \qquad
  t_{m+1}-t_m\leq \Delta T_{\rm filt}.
\)
The filtered state at time \(t_m\) is denoted by
\(
  \psi_m^{\rm filt}.
\)
Initially,
$$
  \psi_0^{\rm filt}
  =
  \psi_0|_{\OmLW}.
$$
On each slab \(I_m=[t_m,t_{m+1}]\), we either train a new neural network
or warm-start from the previous slab.

The initial matching loss is
\begin{equation}\label{eq:sch-initial-loss}
  \mathcal L_{\rm init}^{(m)}(\theta)
  =
  \frac{1}{N_i}
  \sum_{r=1}^{N_i}
  \left|
  \psi_\theta(t_m,y_r)-\psi_m^{\rm filt}(y_r)
  \right|^2.
\end{equation}

The points \(y_r\) can be FFT grid points, random points, or
interpolation points from the filtered spectral grid.

The basic slab loss is
\begin{equation}\label{eq:sch-slab-loss}
  \mathcal L_{\rm Sch}^{(m)}(\theta)
  =
  \mathcal L_{\rm PDE}^{(m)}(\theta)
  +
  \lambda_i\mathcal L_{\rm init}^{(m)}(\theta)
  +
  \lambda_{\rm per}\mathcal L_{\rm per}^{(m)}(\theta)
  +
  \lambda_{\rm reg}\mathcal L_{\rm reg}^{(m)}(\theta).
\end{equation}

A periodic term \(\mathcal L_{\rm per}^{(m)}\) is optional and is
imposed only on the outer boundary of \(\OmLW\) to stabilize the
FFT-based implementation.  

\begin{algorithm}[tbp]
\caption{Phase-space-filtered PINN for the linear Schr\"odinger equation with potential}
\label{alg:sch-pinn-tdpsf}
\begin{algorithmic}[1]
\Require Potential \(V\), initial data \(\psi_0\), final time \(T\),
physical box \(\OmL\), extended box \(\OmLW\), filter parameters from
Algorithm~\ref{alg:sch-precompute}, neural architecture
\(\psi_\theta=a_\theta+\ii b_\theta\), and loss weights.
\Ensure A filtered neural approximation \(\psi_{\theta_m}\) on each time
slab and filtered grid states \(\psi_m^{\rm filt}\).
\State Set $$\psi_0^{\rm filt}=\psi_0$$ on the grid of $$\OmLW.$$
\State Choose time nodes
\[
  0=t_0<t_1<\cdots<t_M=T
\]
with
\[
  t_{m+1}-t_m\leq\Delta T_{\rm filt}.
\]
\For{\(m=0,1,\ldots,M-1\)}
  \State Define the slab \(I_m=[t_m,t_{m+1}]\).
  \State Generate residual collocation points in \(I_m\times\OmLW\),
  with additional samples in \(I_m\times\Bw\) and near the support of
  \(V\).
  \State Generate initial matching points \(\{y_r\}\subset\OmLW\) and
  evaluate or interpolate \(\psi_m^{\rm filt}(y_r)\).
  \State Train \(\theta_m\) by minimizing
  $$\mathcal L_{\rm Sch}^{(m)}$$ in \eqref{eq:sch-slab-loss}.
  \State Evaluate the trained network at the slab endpoint on the FFT
  grid:
  \[
    \psi_{m+1}^{\rm raw}(x_i)
    =
    \psi_{\theta_m}(t_{m+1},x_i).
  \]
  \State Apply the outgoing filter of Algorithm~\ref{alg:sch-filter}:
  \[
    \psi_{m+1}^{\rm filt}
    =
    \mathsf F_{\rm Sch}\psi_{m+1}^{\rm raw}.
  \]
  \State Store \(\psi_{\theta_m}\) as the continuous neural
  representation on \(I_m\times\OmL\).
  \State Use $$\psi_{m+1}^{\rm filt}$$ as the initial condition for the
  next slab.
\EndFor
\end{algorithmic}
\end{algorithm}

\section{Nonlinear Schr\"odinger equations and the TDPSF mechanism}
\label{sec:nls}

We now incorporate nonlinear Schr\"odinger equations.  The model problem is
\begin{equation}\label{eq:nls-general}
  \ii\partial_t\psi
  =
  -\frac12\Delta\psi
  +V(x)\psi
  +\beta |\psi|^{2\sigma}\psi
  +f(t,x),
  \qquad x\in\mathbb R^d,
\end{equation}
where \(\beta\in\mathbb R\), \(\sigma>0\), and \(V\) is real-valued.
The cubic NLS corresponds to \(\sigma=1\).  More generally, the
nonlinear term may be replaced by
\begin{equation}\label{eq:nls-general-g}
  g(t,x,\psi)\psi,
\end{equation}
provided that the far-field behavior in the buffer is asymptotically
linear.

The TDPSF idea is especially natural for NLS.  The solution may be
strongly nonlinear in the physical region \(\OmL\), while the outgoing
radiation near the artificial boundary is often well described by
linear Schr\"odinger wave packets.  Thus the nonlinear equation is
solved in the extended box \(\OmLW\), but the open-boundary treatment is
performed by a phase-space criterion based on outgoing linear radiation
in the buffer

$$
  \Bw=\OmLW\setminus\OmL.
$$

\begin{assumption} 
\label{ass:nls-far-field}
The filter is applied at times when the potential is constant or
negligible in the buffer \(\Bw\), and the nonlinear density in the
buffer is small in the sense that
\begin{equation}\label{eq:nls-buffer-nonlinearity-small}
  \int_{\Bw}|\psi(t,x)|^{2\sigma+2}\,\mathrm dx
\end{equation}
is small at filtering times.  In this regime the outgoing group velocity
is computed from the far-field linear Hamiltonian
\[
  H_\infty(k)=\frac12|k|^2+V_\infty,
  \qquad
  v_g(k)=\nabla_k H_\infty(k)=k.
\]
If \eqref{eq:nls-buffer-nonlinearity-small} is not small, then the
filter should report a warning rather than silently delete the wave
packet.
\end{assumption}

\subsection{PINN residual for NLS}
\label{subsec:nls-pinn-residual}

We write the neural approximation as
\[
  \psi_\theta=a_\theta+\ii b_\theta,
  \qquad
  \rho_\theta=a_\theta^2+b_\theta^2.
\]
For \eqref{eq:nls-general}, define the complex residual
\begin{equation}\label{eq:nls-complex-residual}
  \mathcal R_\theta^{\rm NLS}
  =
  \ii\partial_t\psi_\theta
  +\frac12\Delta\psi_\theta
  -V\psi_\theta
  -\beta |\psi_\theta|^{2\sigma}\psi_\theta
  -f.
\end{equation}
If $$f=f_R+\ii f_I,$$ then the real and imaginary residuals are
\begin{align}
  \mathcal R_{\theta,{\rm NLS}}^R
  &=
  -\partial_t b_\theta
  +\frac12\Delta a_\theta
  -Va_\theta
  -\beta \rho_\theta^{\sigma}a_\theta
  -f_R,
  \label{eq:nls-residual-real}\\
  \mathcal R_{\theta,{\rm NLS}}^I
  &=
  \partial_t a_\theta
  +\frac12\Delta b_\theta
  -Vb_\theta
  -\beta \rho_\theta^{\sigma}b_\theta
  -f_I.
  \label{eq:nls-residual-imag}
\end{align}
On the slab
\(
  I_m=[t_m,t_{m+1}],
\)
the residual loss is approximated by collocation:
\begin{equation}\label{eq:nls-residual-loss}
  \mathcal L_{\rm PDE,NLS}^{(m)}(\theta)
  =
  \frac{1}{N_r}
  \sum_{r=1}^{N_r}
  \left(
  |\mathcal R_{\theta,{\rm NLS}}^R(\tau_r,x_r)|^2
  +
  |\mathcal R_{\theta,{\rm NLS}}^I(\tau_r,x_r)|^2
  \right),
\end{equation}
where
\(
  (\tau_r,x_r)\in I_m\times\OmLW.
\)
The collocation points should be concentrated in the physical box
\(\OmL\), near the support of the potential \(V\), near regions where
\(|\psi_\theta|\) is large, and in the side buffers
\(
  \mathcal B_{w,j}^{s},
  j=1,\ldots,d,\quad s\in\{+,-\}.
\)

The initial matching loss is now
\begin{equation}\label{eq:nls-initial-loss}
  \mathcal L_{\rm init,NLS}^{(m)}(\theta)
  =
  \frac{1}{N_i}
  \sum_{r=1}^{N_i}
  \left|
  \psi_\theta(t_m,y_r)-\psi_m^{\rm filt}(y_r)
  \right|^2.
\end{equation}
Here \(\psi_m^{\rm filt}\) is the filtered state inherited from the
previous slab, and \(y_r\in\OmLW\) are grid points, random points, or
interpolation points from the filtered spectral grid.

The basic NLS slab loss is then
\begin{equation}\label{eq:nls-slab-loss}
  \mathcal L_{\rm NLS}^{(m)}(\theta)
  =
  \mathcal L_{\rm PDE,NLS}^{(m)}(\theta)
  +
  \lambda_i\mathcal L_{\rm init,NLS}^{(m)}(\theta)
  +
  \lambda_{\rm per}\mathcal L_{\rm per}^{(m)}(\theta)
  +
  \lambda_{\rm reg}\mathcal L_{\rm reg}^{(m)}(\theta).
\end{equation}

The term \(\mathcal L_{\rm per}^{(m)}\) is optional and is used only to
make the extended box compatible with FFT-based filtering.  It should
be weighted mildly, because the physical boundary condition is not
periodicity on \(\OmL\).  The role of the phase-space filter is to
remove outgoing radiation in the buffer before it reaches the outer
boundary of \(\OmLW\).

\subsection{Coherent-state outgoing classification}
\label{subsec:nls-coherent-state-classification}

Let \(\Lambda_q\subset\OmLW\) be a position lattice and let
\(
  \Lambda_p\subset[-K_{\max},K_{\max}]^d
\)
be a frequency lattice.  At a filtering time, given a raw endpoint state
\(u\), usually

$$
  u=\psi_{m+1}^{\rm raw},
$$
we compute the coherent-state coefficients
\begin{equation}\label{eq:nls-coefficients}
  c_{q,p}
  =
  (\mathcal G_\sigma u)(q,p),
  \qquad
  (q,p)\in\Lambda_q\times\Lambda_p.
\end{equation}

Equivalently, if \(\phi_{q,p}^{\sigma}\) denotes the coherent state with
center \(q\), momentum \(p\), and width \(\sigma\), then
$$c_{q,p}=\langle u,\phi_{q,p}^{\sigma}\rangle.$$  To classify
\(c_{q,p}\), we first determine whether \(q\) lies in the buffer
\(\Bw\): if \(q\in\mathcal B_{w,j}^{s}\), we use the outward normal
\(n(q)=n_j^s\).  For the standard far-field Schr\"odinger operator,
\(v_g(p)=p\), while for an anisotropic effective mass tensor \(A\), one
uses \(v_g(p)=Ap\).

Given a glancing buffer \(\gamma>0\), the coefficient is classified as displayed below

\[
\begin{array}{lll}
  \text{outgoing}
  &\text{if}&
  q\in\Bw \ \text{and}\ v_g(p)\cdot n(q)\geq \gamma,
  \\[1mm]
  \text{incoming}
  &\text{if}&
  q\in\Bw \ \text{and}\ v_g(p)\cdot n(q)\leq -\gamma,
  \\[1mm]
  \text{glancing or ambiguous}
  &\text{if}&
  q\in\Bw \ \text{and}\ |v_g(p)\cdot n(q)|<\gamma,
  \\[1mm]
  \text{interior}
  &\text{if}&
  q\notin\Bw.
\end{array}
\]

A smooth version replaces the hard outgoing indicator by
\begin{equation}\label{eq:nls-smooth-outgoing-mask}
  \chi_{\rm out}(q,p)
  =
  \eta_{\rm bdry}(q)
  S_\alpha\!\left(v_g(p)\cdot n(q)-\gamma\right),
\end{equation}
where \(\eta_{\rm bdry}\) is a smooth cutoff supported in \(\Bw\), and

$$
  S_\alpha(r)=\frac{1}{1+\exp(-r/\alpha)}
$$

is a smoothed Heaviside function.  The filtered coherent-state
coefficients are then
\begin{equation}\label{eq:nls-filtered-coefficients}
  c_{q,p}^{\rm filt}
  =
  \big(1-\chi_{\rm out}(q,p)\big)c_{q,p}.
\end{equation}

If the coherent states form a tight frame with frame constant \(A_0\),
the reconstructed filtered state is then
\begin{equation}\label{eq:nls-coherent-reconstruction}
  u^{\rm filt}(x)
  =
  A_0^{-1}
  \sum_{(q,p)\in\Lambda_q\times\Lambda_p}
  c_{q,p}^{\rm filt}\phi_{q,p}^{\sigma}(x).
\end{equation}
For a non-tight frame, one uses the dual frame or solves the associated
least-squares reconstruction problem.

\subsection{Original-style NLS TDPSF step}
\label{subsec:nls-original-tdpsf}

The following algorithm describes the NLS phase-space-filter step in the
spirit of the time-dependent phase-space filter.

The interior
propagator may be a spectral method, a finite difference method, or the
slabwise PINN.

\begin{algorithm}[tbp]
\caption{Original-style TDPSF step for nonlinear Schr\"odinger equations}
\label{alg:nls-original-tdpsf}
\begin{algorithmic}[1]
\Require Current filtered state \(u_m^{\rm filt}\) at time \(t_m\),
time interval \([t_m,t_{m+1}]\), extended box \(\OmLW\),
coherent-state lattices \(\Lambda_q,\Lambda_p\), window width
\(\sigma\), glancing buffer \(\gamma\), and tolerance parameters.
\Ensure Filtered state \(u_{m+1}^{\rm filt}\) and diagnostics.
\State Advance the NLS \eqref{eq:nls-general} on \(\OmLW\) from
\(t_m\) to \(t_{m+1}\) using the chosen interior solver.  Denote the
raw endpoint by \(u_{m+1}^{\rm raw}\).
\State Compute the coherent-state coefficients
\[
  c_{q,p}
  =
  (\mathcal G_\sigma u_{m+1}^{\rm raw})(q,p),
  \qquad
  (q,p)\in\Lambda_q\times\Lambda_p.
\]
\State For each \((q,p)\), determine whether \(q\in\Bw\).  If
\(q\notin\Bw\), mark the coefficient as interior and keep it.
\State If \(q\in\Bw\), identify a side buffer \(\mathcal B_{w,j}^{s}\)
containing \(q\), set \(n(q)=n_j^s\), and compute the far-field group
velocity $$v_g(p)=p$$ for the standard Schr\"odinger operator, or
$$v_g(p)=Ap$$ for anisotropic mass tensor \(A\).
\State Classify the coefficient:
\[
\begin{array}{lll}
\text{outgoing} &\text{if}& v_g(p)\cdot n(q)\geq \gamma,\\
\text{incoming} &\text{if}& v_g(p)\cdot n(q)\leq -\gamma,\\
\text{ambiguous/glancing} &\text{if}& |v_g(p)\cdot n(q)|<\gamma.
\end{array}
\]
\State Delete outgoing coefficients by setting
\(c_{q,p}^{\rm filt}=0\).  Keep all other coefficients by setting
$$c_{q,p}^{\rm filt}=c_{q,p}.$$
\State Reconstruct the filtered state using the frame reconstruction
formula
\[
  u_{m+1}^{\rm filt}
  =
  A_0^{-1}
  \sum_{(q,p)\in\Lambda_q\times\Lambda_p}
  c_{q,p}^{\rm filt}\phi_{q,p}^{\sigma}
\]
for a tight frame, or use the corresponding dual-frame reconstruction.
\State Compute the removed coefficient mass and glancing coefficient
mass:
\[
  M_{\rm rem}
  =
  \sum_{(q,p)\in\mathcal O}|c_{q,p}|^2,
  \qquad
  M_{\rm gl}
  =
  \sum_{(q,p)\in\mathcal G}|c_{q,p}|^2,
\]
where \(\mathcal O\) is the outgoing index set and \(\mathcal G\) is the
glancing/ambiguous index set.
\State Compute the nonlinear buffer diagnostic
\[
  N_{\rm buf}
  =
  \int_{\Bw}|u_{m+1}^{\rm raw}(x)|^{2\sigma+2}\,\mathrm dx.
\]
\State If \(M_{\rm gl}\) or \(N_{\rm buf}\) exceeds the prescribed
tolerance, flag the filtering step as unreliable and reduce the time
step, enlarge the buffer, or refine the phase-space lattice.
\State Return $$u_{m+1}^{\rm filt}$$ and the diagnostics.
\end{algorithmic}
\end{algorithm}

\subsection{FFT-window version of the NLS filter}
\label{subsec:nls-fft-window-filter}

For high-dimensional implementations, one may replace the coherent-state
deletion step by the localized FFT-window filter from
Section~\ref{subsec:fft-window-filter}.  For the scalar NLS, define

\[
  P_j^s(k)
  =
  S_\alpha\!\left(v_g(k)\cdot n_j^s-\gamma\right),
  \qquad
  v_g(k)=k
\]
for the standard far-field Schr\"odinger operator.  For anisotropic
mass tensor \(A\), use \(v_g(k)=Ak\).  The outgoing component through side \((j,s)\) is
$$\mathcal O_j^s u=\eta_j^s(x)\mathcal F^{-1}
[P_j^s(k)\mathcal F(\eta_j^s u)(k)](x)$$ and the corresponding
one-side filter is $$\mathsf F_j^s u=u-\mathcal O_j^s u.$$

The complete NLS FFT-window filter is then

\begin{equation}\label{eq:nls-fft-complete-filter}
  \mathsf F_{\rm NLS}u
  =
  \left(
  \prod_{j=1}^{d}
  \mathsf F_j^{-}\mathsf F_j^{+}
  \right)u.
\end{equation}

Equivalently, in additive schematic form, we write

\[
  \mathsf F_{\rm NLS}u
  \approx
  u-
  \sum_{j=1}^d
  \sum_{s\in\{+,-\}}
  \eta_j^s(x)\mathcal F^{-1}
  \left[
    P_j^s(k)\mathcal F(\eta_j^s u)(k)
  \right].
\]

\subsection{NLS TDPSF--PINN algorithm}
\label{subsec:nls-pinn-tdpsf}

We now replace the generic interior propagator in
Algorithm~\ref{alg:nls-original-tdpsf} by a slabwise PINN.  The filtered
state at the beginning of the slab \(I_m=[t_m,t_{m+1}]\) is
\(\psi_m^{\rm filt}\).  The network is trained by minimizing
\eqref{eq:nls-slab-loss}.  Then its endpoint is filtered.

\begin{algorithm}[tbp]
\caption{Phase-space-filtered PINN for nonlinear Schr\"odinger equations}
\label{alg:nls-pinn-tdpsf}
\begin{algorithmic}[1]
\Require Nonlinearity parameters \(\beta,\sigma\), potential \(V\),
initial data \(\psi_0\), final time \(T\), boxes
\(\OmL\subset\OmLW\), coherent-state or FFT-window filter parameters,
neural architecture $$\psi_\theta=a_\theta+\ii b_\theta$$ and loss
weights.
\Ensure Filtered neural approximations \(\psi_{\theta_m}\) on time
slabs and filtered grid states \(\psi_m^{\rm filt}\).
\State Initialize \(\psi_0^{\rm filt}=\psi_0\) on the grid of \(\OmLW\).
\State Choose time nodes
\[
  0=t_0<t_1<\cdots<t_M=T
\]
so that outgoing resolved packets cannot cross the buffer during one
slab.  For the standard far-field Schr\"odinger operator, a typical
choice is
\[
  t_{m+1}-t_m
  \leq
  \Delta T_{\rm filt}
  \leq
  \frac{w}{3K_{\max}}.
\]
For anisotropic mass tensor \(A\), replace \(K_{\max}\) by
\[
  v_{\max}=\sup_{|k|\leq K_{\max}} |Ak|.
\]
\For{\(m=0,1,\ldots,M-1\)}
  \State Define the time slab \(I_m=[t_m,t_{m+1}]\).
  \State Generate residual collocation points in \(I_m\times\OmLW\).
  Oversample the nonlinear interaction region, the potential region,
  and the buffer \(\Bw\).
  \State Generate initial matching points \(\{y_r\}\subset\OmLW\) and
  interpolate \(\psi_m^{\rm filt}(y_r)\).
  \State Train \(\theta_m\) by minimizing the NLS loss
  \(\mathcal L_{\rm NLS}^{(m)}\) in \eqref{eq:nls-slab-loss}.
  \State Evaluate the raw endpoint on the FFT grid:
  \[
    \psi_{m+1}^{\rm raw}(x_i)
    =
    \psi_{\theta_m}(t_{m+1},x_i).
  \]
  \State Apply either the coherent-state filter of
  Algorithm~\ref{alg:nls-original-tdpsf} or the FFT-window filter
  \eqref{eq:nls-fft-complete-filter}:
  \[
    \psi_{m+1}^{\rm filt}
    =
    \mathsf F_{\rm NLS}\psi_{m+1}^{\rm raw}.
  \]
  \State Compute and store the removed coefficient mass, glancing
  energy, high-frequency tail, and nonlinear buffer diagnostic.
  \State If a diagnostic exceeds tolerance, refine the slab: reduce
  \(t_{m+1}-t_m\), enlarge \(w\), refine
  \(\Lambda_q\times\Lambda_p\), or retrain with more buffer collocation
  points.
  \State Store \(\psi_{\theta_m}\) as the neural representation on
  \(I_m\times\OmL\).
  \State Use $$\psi_{m+1}^{\rm filt}$$ as the initial state for the next
  slab.
\EndFor
\end{algorithmic}
\end{algorithm}

\subsection{Practical NLS failure indicators}
\label{subsec:nls-failure-indicators}

The NLS TDPSF--PINN algorithm should report the following warnings.

\begin{enumerate}[label=\textup{(\roman*)}]
\item \emph{Large nonlinear buffer energy.}  If $$\int_{\mathcal B_w}|\psi|^{2\sigma+2}\,\mathrm dx$$ is large, the
far-field linear outgoing classifier is questionable.

\item \emph{Large glancing energy.}  If many coherent states satisfy $$|v_g(p)\cdot n(q)|<\gamma$$ then the
algorithm cannot confidently decide whether they are outgoing or
incoming.

\item \emph{Large high-frequency tail.}  Since \(v_g(k)=k\) for the standard Schr\"odinger equation, high
frequencies may rapidly cross the buffer.

Hence the time-step
restriction must be based on the resolved maximum frequency
\(K_{\max}\), or, in the anisotropic case, on
$$v_{\max}=\sup_{|k|\leq K_{\max}}|v_g(k)|.$$

\item \emph{Large coefficient mass near corners.}  Near corners, a
packet may be outgoing relative to more than one side
\((j,s)\).

A partition of unity or a corner-specific normal should be
used to avoid double counting.

\item \emph{Large mismatch after filtering.}  If \(\|\psi_{m+1}^{\rm raw}-\psi_{m+1}^{\rm filt}\|_{L^2(\Omega_L)}\)
is not small, then the filter is affecting the physical box, indicating
that the buffer is too narrow, the side window is too wide, or the slab
length is too large.
\end{enumerate}

These diagnostics are part of the value of the method: rather than
silently producing a reflected or wrapped solution, the algorithm
identifies when the phase-space classification is unreliable.

\section{Anisotropic wave equations on unbounded domains}
\label{sec:anisotropic}

We next consider linear anisotropic wave systems.  The far-field model
is a constant-coefficient first-order system
\begin{equation}\label{eq:anisotropic-system}
  \partial_t U(t,x)=\mathcal H U(t,x),
  \qquad
  \mathcal H=\sum_{r=1}^d A_r\partial_{x_r}+B,
\end{equation}
where
$$
  U:\mathbb R^d\to\mathbb C^q.
$$

We assume that the matrices \(A_r\in\mathbb C^{q\times q}\) are
Hermitian and that \(B\in\mathbb C^{q\times q}\) is skew-Hermitian.
More generally, a symmetrizable hyperbolic system
\[
  T\partial_t U
  =
  \sum_{r=1}^d A_r\partial_{x_r}U+BU,
  \qquad
  T=T^*>0,
\]

can be reduced to the form \eqref{eq:anisotropic-system} by the change
of variables $$U\mapsto T^{1/2}.$$

The actual PDE may contain scatterers, lower-order terms,
nonlinearities, or variable coefficients inside a compact region
\(B_{\rm int}\Subset\Omega_L\).  The filter, however, is constructed
from the constant-coefficient far-field operator \(\mathcal H\) in the
buffer \(\mathcal B_w=\Omega_{L,w}\setminus\Omega_L\).  Thus the
open-boundary mechanism is determined by the far-field dispersion
relation and by the group velocities of the outgoing branches.

\subsection{Dispersion branches and group velocities}
\label{subsec:aniso-branches-group-velocities}

For a plane wave

$$ (U(t,x)=d(k)\exp\!\big(\ii(k\cdot x-\omega t)\big),
$$
the Fourier symbol of \(\mathcal H\) is
\[
  \widehat{\mathcal H}(k)
  =
  \ii\sum_{r=1}^d k_rA_r+B.
\]
Since \(\widehat{\mathcal H}(k)\) is skew-Hermitian, the matrix

\begin{equation}\label{eq:M-k-Hermitian}
  M(k)=-\ii\widehat{\mathcal H}(k)
\end{equation}
is Hermitian.  Let
\begin{equation}\label{eq:dispersion-eigenproblem}
  M(k)d_\ell(k)=\omega_\ell(k)d_\ell(k),
  \qquad
  \ell=1,\ldots,q,
\end{equation}
where the eigenvectors \(d_\ell(k)\) are chosen orthonormally whenever
the eigenvalues are simple.  The group velocity on branch \(\ell\) is

\begin{equation}\label{eq:anisotropic-group-velocity}
  v_{g,\ell}(k)=\nabla_k\omega_\ell(k).
\end{equation}
If the eigenvalue \(\omega_\ell(k)\) is simple, then the
Hellmann--Feynman formula gives
\begin{equation}\label{eq:hellmann-feynman}
  \partial_{k_r}\omega_\ell(k)
  =
  d_\ell(k)^*\partial_{k_r}M(k)d_\ell(k),
  \qquad r=1,\ldots,d.
\end{equation}

At repeated eigenvalues, or near crossings, one should use spectral
cluster projectors rather than individual eigenvectors.

\begin{definition}
\label{def:aniso-outgoing-branch}
For the side \((j,s)\), with \(j\in\{1,\ldots,d\}\) and
\(s\in\{+,-\}\), the branch \(\ell\) is outgoing at frequency \(k\) if
\begin{equation}\label{eq:anisotropic-outgoing-condition}
  v_{g,\ell}(k)\cdot n_j^s>0.
\end{equation}

\end{definition}

This group-velocity criterion distinguishes isotropic from anisotropic
propagation: in isotropic models, phase and group velocities are aligned,
whereas in anisotropic models they need not be.  Therefore filtering by
the sign of \(k\cdot n_j^s\) is generally incorrect, then we must filter by
the sign of \(v_{g,\ell}(k)\cdot n_j^s\).

\subsection{Branchwise outgoing projectors}
\label{subsec:aniso-frequency-projectors}

Let
\(
  \Pi_\ell(k)=d_\ell(k)d_\ell(k)^*
\)
be the spectral projector onto branch \(\ell\) when the eigenvalue is
simple.  Near multiple eigenvalues, \(\Pi_\ell(k)\) should be replaced
by the projector onto the corresponding spectral cluster.

Choose a glancing buffer \(\gamma>0\) and a smoothing width
\(\alpha>0\).  For branch \(\ell\), side direction \(j\), and sign
\(s\), define the scalar branch mask
\begin{equation}\label{eq:anisotropic-branch-mask}
  p_{\ell,j}^{s}(k)
  =
  S_\alpha\!\left(v_{g,\ell}(k)\cdot n_j^s-\gamma\right),
\end{equation}
where

$$
 S_\alpha(r)=\frac{1}{1+\exp(-r/\alpha)}.
$$

The branchwise outgoing Fourier multiplier is now
\begin{equation}\label{eq:anisotropic-branch-projector}
  P_{\ell,j}^{s}(k)
  =
  p_{\ell,j}^{s}(k)\Pi_\ell(k).
\end{equation}
The total outgoing projector for side \((j,s)\) is then

\begin{equation}\label{eq:anisotropic-frequency-projector}
  P_j^s(k)
  =
  \sum_{\ell=1}^{q}P_{\ell,j}^{s}(k)
  =
  \sum_{\ell=1}^{q}
  p_{\ell,j}^{s}(k)\Pi_\ell(k).
\end{equation}
Equivalently, if \(D(k)\) is the unitary matrix whose columns are the
eigenvectors \(d_\ell(k)\), then
\[
  P_j^s(k)
  =
  D(k)
  \operatorname{diag}
  \big(
  p_{1,j}^{s}(k),\ldots,p_{q,j}^{s}(k)
  \big)
  D(k)^*.
\]
The precise placement of \(D(k)\) and \(D(k)^*\) depends only on whether
the eigenvectors are stored as columns or rows. Throughout the work,
\(P_j^s(k)\) denotes the projector acting on the physical vector
\(\widehat U(k)\in\mathbb C^q\).

\subsection{Outgoing wave filter for anisotropic systems}
\label{subsec:aniso-outgoing-filter}

Let \(U:\OmLW\to\mathbb C^q\) be a grid function.  For side \((j,s)\),
let \(\eta_j^s\) be the side window supported in
\(\mathcal B_{w,j}^{s}\).  If branch-dependent localization is needed,
we write \(\eta_{\ell,j}^{s}\) and
\(\mathcal B_{w,\ell,j}^{s}\).  Otherwise the spatial window is shared
by all branches and we set

$$
  \eta_{\ell,j}^{s}=\eta_j^s.
$$

Using the total outgoing projector \(P_j^s(k)\), the outgoing component
through side \((j,s)\) is
\begin{equation}\label{eq:anisotropic-outgoing-component}
  \mathcal O_j^s U
  =
  \eta_j^s(x)\mathcal F^{-1}
  \left[
    P_j^s(k)\mathcal F(\eta_j^s U)(k)
  \right](x),
\end{equation}
in which the Fourier transform is applied componentwise to the vector field
\(U\), and \(P_j^s(k)\) acts as a \(q\times q\) matrix multiplier on the
Fourier coefficient \(\widehat U(k)\).

The corresponding one-side filter is

\begin{equation}\label{eq:anisotropic-one-side-filter}
  \mathsf F_j^s U
  =
  U-\mathcal O_j^s U.
\end{equation}
The complete anisotropic phase-space filter is as follows
\begin{equation}\label{eq:anisotropic-complete-filter}
  \mathsf F_{\rm aniso}U
  =
  \left(
  \prod_{j=1}^d
  \mathsf F_j^-\mathsf F_j^+
  \right)U.
\end{equation}

Equivalently, in additive schematic form, we find

\[
  \mathsf F_{\rm aniso}U
  \approx
  U-
  \sum_{j=1}^d
  \sum_{s\in\{+,-\}}
  \eta_j^s(x)\mathcal F^{-1}
  \left[
    P_j^s(k)\mathcal F(\eta_j^s U)(k)
  \right].
\]
If branch-dependent windows are used, this may be written as
\[
  \mathsf F_{\rm aniso}U
  \approx
  U-
  \sum_{j=1}^d
  \sum_{s\in\{+,-\}}
  \sum_{\ell=1}^q
  \eta_{\ell,j}^{s}(x)\mathcal F^{-1}
  \left[
    P_{\ell,j}^{s}(k)
    \mathcal F(\eta_{\ell,j}^{s} U)(k)
  \right].
\]

\begin{algorithm}[tbp]
\caption{Precomputation for anisotropic phase-space filters}
\label{alg:aniso-precompute}
\begin{algorithmic}[1]
\Require Far-field matrices \(A_1,\ldots,A_d,B\), physical half-width
\(L\), buffer width \(w\), tolerance \(\delta\), glancing buffer
\(\gamma\), mask smoothing width \(\alpha\), FFT grid, and maximum
resolved frequency \(K_{\max}\).
\Ensure Side windows \(\eta_j^s\), outgoing projectors $$P_j^s(k_r),$$
filter time step \(\Delta T_{\rm filt}\), and ambiguity diagnostics.
\State Build the extended grid on
\[
  \OmLW=[-L-w,L+w]^d
\]
and the corresponding Fourier grid \(\{k_r\}\).
\State Construct the side windows \(\eta_j^s\) as in
Section~\ref{subsec:side-window-construction}, for all
\(j=1,\ldots,d\) and \(s\in\{+,-\}\).
\For{each Fourier grid point \(k_r\)}
  \State Form the Hermitian matrix
  \[
    M(k_r)
    =
    -\ii
    \left(
    \ii\sum_{m=1}^d k_{r,m}A_m+B
    \right).
  \]
  \State Compute the eigenpairs
  \[
    M(k_r)d_\ell(k_r)
    =
    \omega_\ell(k_r)d_\ell(k_r),
    \qquad
    \ell=1,\ldots,q.
  \]
  \State Compute group velocities
  \[
    v_{g,\ell}(k_r)=\nabla_k\omega_\ell(k_r).
  \]
  For simple branches, use the Hellmann--Feynman formula
  \eqref{eq:hellmann-feynman}; near repeated eigenvalues, use finite
  differences or spectral-cluster projectors.
  \For{each side \((j,s)\)}
    \State Compute branch masks

    \[
      p_{\ell,j}^{s}(k_r)
      =
      S_\alpha\!\left(
      v_{g,\ell}(k_r)\cdot n_j^s-\gamma
      \right),
      \qquad
      \ell=1,\ldots,q.
    \]
    \State Assemble the total outgoing projector

    \[
      P_j^s(k_r)
      =
      \sum_{\ell=1}^q
      p_{\ell,j}^{s}(k_r)
      d_\ell(k_r)d_\ell(k_r)^*.
    \]
    \State Store the glancing indicator
    \[
      g_{\ell,j}^{s}(k_r)
      =
      \mathbf 1_{\{
      |v_{g,\ell}(k_r)\cdot n_j^s|\leq \gamma
      \}}.
    \]
  \EndFor
\EndFor
\State Set

\[
  v_{\max}
  =
  \sup_{|k_r|\leq K_{\max}}
  \max_{\ell=1,\ldots,q}
  |v_{g,\ell}(k_r)|.
\]
\State Choose

\[
  \Delta T_{\rm filt}\leq \frac{w}{3v_{\max}}.
\]
\State Verify the uncertainty condition

\[
  C_1 k_b^{-1}\sqrt{\log(\delta^{-1})}
  \leq
  \sigma
  \leq
  C_2\frac{w}{\sqrt{\log(\delta^{-1})}}.
\]
If it fails, adjust \(w\), \(\gamma\), or \(\sigma\).
\end{algorithmic}
\end{algorithm}

\begin{algorithm}[tbp]
\caption{One application of the anisotropic outgoing wave filter}
\label{alg:aniso-filter}
\begin{algorithmic}[1]
\Require Grid values \(U(x_i)\in\mathbb C^q\), side windows
$$\eta_j^s,$$ and matrix projectors \(P_j^s(k_r)\).
\Ensure Filtered grid values \(U^{\rm filt}\).
\State Set \(U^{(0)}=U\) and \(r=0\).
\For{\(j=1,\ldots,d\)}
  \For{\(s\in\{+,-\}\)}
    \State Localize each component:
    \[
      Z(x_i)=\eta_j^s(x_i)U^{(r)}(x_i).
    \]
    \State Compute the componentwise FFT
    \[
      \widehat Z(k_r)=\mathcal F Z(k_r)\in\mathbb C^q.
    \]
    \For{each Fourier grid point \(k_r\)}
      \State Project onto outgoing branches:
      \[
        \widehat Z_{\rm out}(k_r)
        =
        P_j^s(k_r)\widehat Z(k_r).
      \]
    \EndFor
    \State Transform back and localize again:
    \[
      O(x_i)
      =
      \eta_j^s(x_i)
      \mathcal F^{-1}\widehat Z_{\rm out}(x_i).
    \]
    \State Remove outgoing components:
    \[
      U^{(r+1)}(x_i)=U^{(r)}(x_i)-O(x_i).
    \]
    \State Set \(r\gets r+1\).
  \EndFor
\EndFor
\State Return $$U^{\rm filt}=U^{(r)}.$$
\end{algorithmic}
\end{algorithm}

\subsection{PINN residual for anisotropic systems}
\label{subsec:aniso-pinn-residual}

Let
\(
  U_\theta:I_m\times\OmLW\to\mathbb C^q
\)

be a neural approximation.  For the far-field model
\eqref{eq:anisotropic-system}, the residual is

\begin{equation}\label{eq:aniso-residual-constant}
  \mathcal R_\theta
  =
  \partial_tU_\theta
  -
  \sum_{r=1}^d A_r\partial_{x_r}U_\theta
  -
  BU_\theta.
\end{equation}
For a more general interior model
\begin{equation}\label{eq:aniso-residual-variable-model}
  \partial_t U
  =
  \mathcal H_{\rm int}(x,\partial_x)U+F(t,x),
\end{equation}
we define

\begin{equation}\label{eq:aniso-residual-variable}
  \mathcal R_\theta
  =
  \partial_tU_\theta
  -
  \mathcal H_{\rm int}(x,\partial_x)U_\theta
  -
  F.
\end{equation}
The important requirement is that \(\mathcal H_{\rm int}\) coincides
with the constant-coefficient far-field operator \(\mathcal H\) in the
buffer where the filter is applied.

The residual loss on the slab \(I_m=[t_m,t_{m+1}]\) is
\begin{equation}\label{eq:aniso-residual-loss}
  \mathcal L_{\rm PDE}^{(m)}(\theta)
  =
  \frac{1}{N_r}
  \sum_{r=1}^{N_r}
  \left|
  \mathcal R_\theta(\tau_r,x_r)
  \right|_{\mathbb C^q}^2,
\end{equation}

where
\(
  (\tau_r,x_r)\in I_m\times\OmLW.
\)

The initial matching loss is

\begin{equation}\label{eq:aniso-initial-loss}
  \mathcal L_{\rm init}^{(m)}(\theta)
  =
  \frac{1}{N_i}
  \sum_{r=1}^{N_i}
  \left|
  U_\theta(t_m,y_r)-U_m^{\rm filt}(y_r)
  \right|_{\mathbb C^q}^2.
\end{equation}
The slab loss is
\begin{equation}\label{eq:aniso-slab-loss}
  \mathcal L_{\rm aniso}^{(m)}(\theta)
  =
  \mathcal L_{\rm PDE}^{(m)}(\theta)
  +
  \lambda_i\mathcal L_{\rm init}^{(m)}(\theta)
  +
  \lambda_{\rm per}\mathcal L_{\rm per}^{(m)}(\theta)
  +
  \lambda_{\rm reg}\mathcal L_{\rm reg}^{(m)}(\theta).
\end{equation}

\begin{algorithm}[tbp]
\caption{Phase-space-filtered PINN for anisotropic wave systems}
\label{alg:aniso-pinn-tdpsf}
\begin{algorithmic}[1]
\Require Interior operator \(\mathcal H_{\rm int}\), far-field matrices
\(A_r,B\), initial data \(U_0\), final time \(T\), boxes
$$\OmL\subset\OmLW,$$ filter parameters from
Algorithm~\ref{alg:aniso-precompute}, neural architecture \(U_\theta\),
and loss weights.
\Ensure Filtered slabwise neural approximations \(U_{\theta_m}\) and
filtered grid states \(U_m^{\rm filt}\).
\State Set \(U_0^{\rm filt}=U_0\) on the grid of \(\OmLW\).
\State Choose time nodes
\[
  0=t_0<t_1<\cdots<t_M=T
\]
with
\[
  t_{m+1}-t_m\leq\Delta T_{\rm filt}.
\]
\For{\(m=0,1,\ldots,M-1\)}
  \State Define the slab \(I_m=[t_m,t_{m+1}]\).
  \State Generate residual collocation points in \(I_m\times\OmLW\),
  with additional points near scatterers and in the buffer \(\Bw\).
  \State Generate initial matching points \(\{y_r\}\subset\OmLW\) and
  evaluate or interpolate $$U_m^{\rm filt}(y_r).$$
  \State Train \(\theta_m\) by minimizing
  \(\mathcal L_{\rm aniso}^{(m)}\) in \eqref{eq:aniso-slab-loss}.
  \State Evaluate the trained network on the FFT grid:
  \[
    U_{m+1}^{\rm raw}(x_i)
    =
    U_{\theta_m}(t_{m+1},x_i).
  \]
  \State Apply the anisotropic outgoing filter of
  Algorithm~\ref{alg:aniso-filter}:
  \[
    U_{m+1}^{\rm filt}
    =
    \mathsf F_{\rm aniso}U_{m+1}^{\rm raw}.
  \]
  \State Compute diagnostics: removed outgoing energy, glancing energy,
  and residual error in \(\OmL\).
  \State Store \(U_{\theta_m}\) as the neural representation on
  \(I_m\times\OmL\).
  \State Use \(U_{m+1}^{\rm filt}\) as the initial condition for the
  next slab.
\EndFor
\end{algorithmic}
\end{algorithm}

\begin{revision}
\emergencystretch=1em
\section{Numerical comparisons and long-time stability}
\label{sec:results}

\subsection{Common comparison protocol}
\label{sec:comparison-protocol}

Within each of the tests, the boundary methods use the same physical problem
and observation region. Their neural comparisons share the field
architecture, collocation rule, initial random seed, and optimizer
budget.

\subsubsection{Long-time stability and accuracy}
\label{sec:numerical-stability}
In the work, stability is assessed through long-time propagation of the prescribed
initial-value problems.

Let $0=t_0<t_1<\cdots<t_N=T$ denote the recorded times.
For each boundary method and training seed, let $u_h(t_n,x)$ denote
the computed PINN field at time $t_n$, including the endpoint filtering
step for TDPSF. Let $u_{\rm ref}(t_n,x)$ denote the independently
computed reference field for the same physical problem. The notation
$u_h(t_n)$ and $u_{\rm ref}(t_n)$ denotes these fields as functions
of the spatial variable $x$.

In Test~1, the compared field is the complex wave function $u=\psi$.
In Tests~2 and~3, it is the acoustic state
$u=(p,v_1,v_2)^{\mathsf T}$, consisting of pressure and the two velocity
components. In Test~3, both compared fields are total fields, obtained
by adding the same tabulated incident field to the corresponding
computed scattered field.

Let $\Omega_{\rm obs}$ denote the observation region specified for
each test. We define spatial error norm as 
$$
\begin{gathered}
\|u_h(t_n)-u_{\rm ref}(t_n)\|_{L^2(\Omega_{\rm obs})}^2
=\vphantom{\int}\phantom{\bigg|}\\[1ex]
\begin{cases}
 \displaystyle\int_{\Omega_{\rm obs}}
 |\psi_h(t_n,x)-\psi_{\rm ref}(t_n,x)|^2\,dx,
 &\text{for Test~1},\\[2ex]
 \begin{aligned}
 &\displaystyle\int_{\Omega_{\rm obs}}
 \Bigl(|p_h(t_n,x)-p_{\rm ref}(t_n,x)|^2\\[-0.3ex]
 &\displaystyle\qquad+
 \sum_{j=1}^{2}|v_{j,h}(t_n,x)-v_{j,\rm ref}(t_n,x)|^2
 \Bigr)\,dx,
 \end{aligned}
 &\text{for Tests~2 and~3}.
\end{cases}
\end{gathered}
$$

Let $\Omega_{\rm phys}$ denote the physical region, excluding the
absorbing buffer. In Test~1, $\Omega_{\rm phys}=\Omega_{\rm obs}$.
Test~2 specifies the two regions separately. We use the fixed scale

$$
S_{\rm ref}=
\begin{cases}
 \displaystyle\|u_{\rm ref}(0)\|_{L^2(\Omega_{\rm phys})},
 &\text{in both Tests~1 and~2},\\[1.5ex]
 \displaystyle\max_{0\leq m\leq N}
 \|u_{\rm ref}(t_m)\|_{L^2(\Omega_{\rm obs})},
 &\text{for Test~3}.
\end{cases}
$$
Accuracy at each recorded endpoint is measured by the following formula
\begin{equation}\label{eq:longtime-error}
 E(t_n)=
 \frac{\|u_h(t_n)-u_{\rm ref}(t_n)\|_{L^2(\Omega_{\rm obs})}}
 {S_{\rm ref}}.
\end{equation}

  In our computations, Automatic Differentiation
computes the loss gradients with respect to the trainable parameters. 

The common spatial transition is the $C^2$ quintic function, defined as follows
\[
\vartheta(z)=
\begin{cases}
 0,&z\leq0,\\
 z^3(10-15z+6z^2),&0<z<1,\\
 1,&z\geq1.
\end{cases}
\]

At a fixed filtering time, let $u^{\rm raw}$ denote the grid array
supplied to the filter. The filtered variable is the complex NLS
envelope in Test~1, the acoustic state $(p,v_1,v_2)^{\mathsf T}$ in
Test~2, and the scattered acoustic state in Test~3. For Test~3,
$u^{\rm raw}$ contains the neural scattered field at grid points in
the annular buffer and zero elsewhere on the Cartesian filtering grid.

Let $J$ be the number of spatial windows, and let $j=1,\ldots,J$
specify their order of application. We use $J=4$ for the four boundary
strips in Tests~1 and~2 and $J=16$ for the overlapping angular windows
in Test~3. The common index $j=1,\ldots,16$ corresponds to the angular
window index $\ell=j-1$ used in the description of Test~3.

For each $j$, the real-valued spatial window $\eta_j(x)$ is supported
in the corresponding part of the buffer, with $0\leq\eta_j(x)\leq1$.
For any grid array $q$, we have
\[
 (\eta_jq)(x_\ell)=\eta_j(x_\ell)q(x_\ell),
\]
where $x_\ell$ is a grid node. Let $n_j$ denote the outward unit normal
associated with window $j$. For an angular window in Test~3, it is the
outward radial unit vector at the window's central angle.

The operator $P_j$ acts on Fourier coefficients and weights wave
components according to their outward group velocity relative to
$n_j$. Writing $\widehat q=\mathcal Fq$, its action is then
\[
 (P_j\widehat q)(k)=P_j(k)\widehat q(k),
\]
where $k$ is a discrete Fourier wave vector. The outgoing normal
group-velocity threshold is $\gamma=0.2$ in all three tests.

In Test~1, filtering acts on the envelope obtained by removing the
carrier $k_0=(-1,4)$. The scalar multiplier is
\begin{equation}\label{W1}
 \begin{aligned}
 P_j(k)&=\left[1+\exp\!\left(
 -\frac{n_j\cdot v_g(k+k_0)-0.2}{0.05}\right)\right]^{-1},\\
 v_g(k)&=(k_1+0.4k_2,\;0.4k_1+k_2).
 \end{aligned}
\end{equation}

In the acoustic tests, for $k\neq0$, we set $e=k/|k|$.
For a Fourier coefficient $\widehat q=(\widehat p,\widehat v)^{\mathsf T}$,
with $\widehat v\in\mathbb C^2$, define
\[
\begin{aligned}
 P_j(k)\begin{pmatrix}\widehat p\\\widehat v\end{pmatrix}
 &=\frac12\sum_{\sigma\in\{-1,1\}}
 \frac{\widehat p+\sigma e\cdot\widehat v}
 {1+\exp[-(\sigma n_j\cdot e-0.2)/\alpha]}
 \begin{pmatrix}1\\\sigma e\end{pmatrix},\\
 \mbox{where}&\\
 \alpha&=
 \begin{cases}
 0.1,&\text{Test~2},\\
 0.05,&\text{Test~3},
 \end{cases}
\end{aligned}
\]
in which $\sigma=\pm1$ labels the two acoustic branches, whose group
velocities are $\sigma e$. 

In Tests~2 and~3, we set $P_j(0)=0$.
In Test~1, the zero envelope frequency represents the physical
carrier $k_0$ and retains the multiplier prescribed by
Equation~\eqref{W1}.

Starting from the unfiltered array, the sequential grid update is
\begin{equation}\label{eq:numerical-sequential-filter}
\begin{aligned}
 u^{(0)}&=u^{\rm raw},\\
 K_j&=\eta_j\mathcal F^*P_j\mathcal F\eta_j,
 \qquad j=1,\ldots,J,\\
 u^{(j)}&=
 \begin{cases}
 (I-K_j)u^{(j-1)},&\text{for Tests~1 and~2},\\
 \operatorname{Re}\!\bigl((I-K_j)u^{(j-1)}\bigr),&\text{in Test~3},
 \end{cases}
 \quad j=1,\ldots,J,\\
 u^{\rm filt}&=
 \begin{cases}
 u^{(J)},&\text{in Tests~1 and~3},\\
 \operatorname{Re}u^{(J)},&\text{Test~2},
 \end{cases}
\end{aligned}
\end{equation}
in which $I$ is the identity on grid arrays, $\mathcal F$ is the unitary
discrete Fourier transform applied componentwise, and $\mathcal F^*$
is its inverse.  

In Test~1, the 4 windows are applied on the left, right, bottom,
and top sides, in that order. In Test~2, the order is right, left,
top, and bottom. In Test~3, the angular windows are applied in the
order $\ell=0,\ldots,15$. Each window acts on the array produced by
the preceding window.

Since $0\leq P_j\leq I$ and $0\leq\eta_j\leq1$, we have
\[
 0\leq K_j\leq I,\qquad \|I-K_j\|_2\leq1,
\]
in which the operator inequalities are in the Hermitian positive-semidefinite
sense, and $\|\cdot\|_2$ is the operator norm induced by the Euclidean
grid norm.  

The logistic transition is centered at outward normal group velocity
$0.2$. At zero normal velocity, its value is
$[1+\exp(0.2/\alpha)]^{-1}$, where $\alpha=0.05$ in Tests~1 and~3
and $\alpha=0.1$ in Test~2.

Test~1 transfers its filtered envelope by a periodic quintic spline.
Test~2 transfers the endpoint grid by cubic interpolation. In Test~3,
TDPSF transfers the complete neural endpoint with the interpolated
grid removal subtracted. Hence, the removal is
$u^{\rm raw}-u^{\rm filt}$ on the filtering grid.

For the loss formulas below, $\operatorname{mse}(R)$ denotes the mean
of squared absolute values over every sample and component of the
array $R$.

\Needspace{14\baselineskip}
\subsection{Test 1:}
\label{sec:results-test1}
This test   compares TDPSF--PINN with PML--PINN over $0\leq t\leq100$.

\subsubsection{The model}
We consider the Nonlinear Schr\"odinger Equation (NLS)
\begin{align}
 \ii\psi_t&=-\tfrac12(\psi_{xx}+0.8\psi_{xy}+\psi_{yy})
                +\beta|\psi|^2\psi,\qquad \beta=-1,\label{eq:test1-model}\\
 \psi(0,x,y)&=0.06\exp\!\left[-\frac{(x-4)^2+(y+72)^2}{72}\right]
                   e^{\ii(-x+4y)}.\label{eq:test1-initial}
\end{align}
The linear dispersion and group velocity are
\begin{equation}\label{eq:test1-dispersion}
 \omega(k)=\tfrac12(k_x^2+0.8k_xk_y+k_y^2),\qquad
 v_g(k)=(k_x+0.4k_y,\,0.4k_x+k_y).
\end{equation}
At the carrier $k_0=(-1,4)$, we have
\[
 \omega_0=6.9,\qquad v_g(k_0)=(0.6,3.6),\qquad v_{g,x}(k_0)=0.6>0.
\]
The normal component of the phase-velocity vector is
\[
 v_{{\rm phase},x}
 =\frac{\omega_0k_{0,x}}{|k_0|^2}
 =-\frac{6.9}{17}<0.
\]

The physical and observation rectangle is
$$\Omega_{\rm obs}=(-24,24)\times(-128,128).$$ A buffer of width $w=24$
gives the computational rectangle $$\Omega_{\rm comp}=(-48,48)\times(-152,152).$$ 
The carrier trajectory $$(4,-72)+t(0.6,3.6)$$ reaches the right boundary
at $t=100/3$ and the top boundary at $t=500/9$, as thus it encounters the right
boundary first. The positive-definite dispersion matrix $$G=\left(\begin{smallmatrix}1&0.4\\0.4&1\end{smallmatrix}\right)$$
has eigenvalues $0.6$ and $1.4$.

We use a complex-coordinate PML for the nonlinear
Schr\"odinger equation~\cite{Zheng2007}.
Let $$\rho=e^{\ii\pi/4}$$ and define
\[
\begin{aligned}
	\sigma_x(x)
	&=\sigma_{\max}
	\left(\frac{(|x|-L_x)_+}{w}\right)^3,
	& L_x&=24,\\
	\sigma_y(y)
	&=\sigma_{\max}
	\left(\frac{(|y|-L_y)_+}{w}\right)^3,
	& L_y&=128,
\end{aligned}
\]
where $(s)_+=\max\{s,0\}$, $w$ is the buffer width,
and $\sigma_{\max}$ is the layer-strength parameter.
The stretching coefficients are as follows
\[
a_x(x)=\frac{1}{1+\rho\sigma_x(x)},
\qquad
a_y(y)=\frac{1}{1+\rho\sigma_y(y)}.
\]
Define the differential operators as 
\[
D_x\psi=a_x(x)\partial_x\psi,
\qquad
D_y\psi=a_y(y)\partial_y\psi.
\]
The layer equation is as follows
\begin{equation}\label{eq:test1-layer}
	\ii\psi_t
	=-\tfrac12
	\left\{D_x^2+D_y^2+0.4(D_xD_y+D_yD_x)\right\}\psi
	-|\psi|^2\psi.
\end{equation}

  The nonlinear term is appended as in the related extension in
\cite[Eq.~(3.8)]{Dohnal2009} and ~\cite{Zheng2007}.
The outer PML closure is homogeneous Dirichlet.
The layer strength is $\sigma_{\max}=0.5$ for all 4 seeds.

 In a right-side layer with $a_y=1$ and constant $a_x=a$, we have the dispersion relation
\begin{equation}\label{eq:test1-frozen-pml}
 \omega_{\rm PML}(k)=\tfrac12(a^2k_x^2+0.8ak_xk_y+k_y^2).
\end{equation}

We write
$$
\psi(t,x,y)
=e^{\ii(-x+4y-6.9t)}A(t,x,y),
$$
where $A$ is the complex envelope represented by the PINN.
For TDPSF--PINN, the envelope equation is
$$
\ii\bigl(A_t+0.6A_x+3.6A_y\bigr)
+\tfrac12\bigl(A_{xx}+0.8A_{xy}+A_{yy}\bigr)
+|A|^2A=0.
$$
TDPSF--PINN minimizes this residual throughout the
computational rectangle, including the buffer.

We apply 4 spatial windows successively.
The index $j=1,\ldots,4$ labels their order of application.
For each window, let $z_j(x,y)$ be the signed distance
from the corresponding physical boundary, positive toward
the buffer, and let $n_j$ be its outward unit normal. We then have the table
$$
\begin{array}{c|c|c|c}
	j & \text{buffer strip} & z_j(x,y) & n_j\\ \hline
	1 & \text{left}   & -x-24  & (-1,0)\\
	2 & \text{right}  &  x-24  & (1,0)\\
	3 & \text{bottom} & -y-128 & (0,-1)\\
	4 & \text{top}    &  y-128 & (0,1)
\end{array}
$$
Using the quintic step $\vartheta$ defined above, we set then
$$
\eta_j(x,y)
=
\vartheta\!\left(\frac{z_j(x,y)}{0.6w}\right)
\vartheta\!\left(\frac{w-z_j(x,y)}{0.1w}\right),
\qquad j=1,\ldots,4.
$$

Let $k=(k_x,k_y)$ denote a discrete Fourier wave vector
of the envelope. We define $P_j$ as in \eqref{W1}.

At each filtering time, let $A^{\rm raw}$ denote the grid
values of the computed PINN envelope before filtering.
Starting from $A^{(0)}=A^{\rm raw}$, we put
$$
A^{(j)}
=
A^{(j-1)}
-
\eta_j\mathcal F^{-1}
\left[
P_j\,\mathcal F\!\left(\eta_j A^{(j-1)}\right)
\right],
\qquad j=1,\ldots,4,
$$
and then set $A^{\rm filt}=A^{(4)}$.

Filtering is performed every
$\Delta t_{\rm filter}=0.5$ on a grid with $128$ points
in the $x$ direction and $512$ points in the $y$ direction.
The maximum resolved normal group speed is approximately
$v_{\max}=10.520$, and we obtain then
\[
v_{\max}\Delta t_{\rm filter}
=5.260<8=\frac{w}{3}.
\]

At the outer boundary, the envelope and its first
coordinate derivatives satisfy the periodic conditions
\[
\begin{aligned}
	A(t,-48,y)&=A(t,48,y),&
	\partial_xA(t,-48,y)&=\partial_xA(t,48,y),\\
	A(t,x,-152)&=A(t,x,152),&
	\partial_yA(t,x,-152)&=\partial_yA(t,x,152).
\end{aligned}
\]

We divide the simulation into time intervals
$[t_n,t_{n+1}]$, called time slabs, where
\[
t_n=n\Delta t,\qquad \Delta t=0.5.
\]
On the $n$th slab, write
\[
t=t_n+\tau,\qquad 0\leq\tau\leq\Delta t.
\]

Let $A_n(x,y)$ denote the periodic quintic spline interpolant
of the stored envelope values at time $t_n$.

For the first slab, the stored values come from the
prescribed initial envelope.
For subsequent slabs, they come from the preceding
computed endpoint, including the TDPSF filtering step
when that method is used. 

The carrier group velocity is
$v_g(k_0)=(0.6,3.6)$.
We represent the envelope on the current slab by
\[
A_\theta(t_n+\tau,x,y)
=
A_n(x-0.6\tau,\;y-3.6\tau)
+\tau N_\theta(\tau,z_x,z_y),
\]
where $\theta$ denotes the collection of trainable real
parameters and $N_\theta$ is a complex-valued network.

The network uses the moving, rescaled spatial coordinates as follows
$$
z_x=\frac{x-4-0.6(t_n+\tau)}{18},
\qquad
z_y=\frac{y+72-3.6(t_n+\tau)}{18}.
$$

The network contains $32$ Gaussian functions, indexed by
$\ell=1,\ldots,32$. Each Gaussian is multiplied by the
six spatial monomials
\[
1,\quad z_x,\quad z_y,\quad z_x^2,\quad z_xz_y,\quad z_y^2.
\]
An explicit expression for the network is
\[
\begin{aligned}
	N_\theta(\tau,z_x,z_y)
	={}&0.06
	\sum_{\ell=1}^{32}
	\sum_{\substack{a,b\in\{0,1,2\}\\a+b\leq2}}
	\bigl(c_{\ell,a,b}^{(0)}
	+\tau c_{\ell,a,b}^{(1)}\bigr)
	z_x^a z_y^b\\
	&\qquad{}\times
	\exp\left[
	-\frac12
	\left\|
	B_\ell
	\begin{pmatrix}\tau\\z_x\\z_y\end{pmatrix}
	+b_\ell
	\right\|_2^2
	\right],
\end{aligned}
\]
in which $B_\ell\in\mathbb R^{2\times3}$ and
$b_\ell\in\mathbb R^2$ are the trainable parameters of
the $\ell$th Gaussian.
For each of the six monomials,
$c_{\ell,a,b}^{(0)},c_{\ell,a,b}^{(1)}\in\mathbb C$
are trainable output coefficients.  

Thus the number of real trainable parameters is computed as
\[
\underbrace{32(2\cdot3+2)}_{\text{Gaussian parameters}}
+
\underbrace{32\cdot6\cdot2\cdot2}_{\text{output coefficients}}
= 1024.
\]

During our optimization, the envelope is represented as
$a_\theta=A_\theta/0.06$.

The physical wave function is reconstructed by the formula
\[
\psi_\theta(t,x,y)
=
e^{\ii(-x+4y-6.9t)}A_\theta(t,x,y).
\]
 
Each time slab has length $\Delta t=0.5$.
On each slab, we draw $1536$ spatial points for evaluating
the PDE residual, as displayed below
\[
\begin{array}{c|l}
	\text{number of points} & \text{sampling rule}\\ \hline
	384 & \text{uniform in }\Omega_{\rm comp},\\
	384 & \text{uniform in }\Omega_{\rm obs},\\
	768 & \text{sampled using the preceding computed envelope}.
\end{array}
\]

For the last group of $768$ points, let
\[
\boldsymbol{x}_\ell=(x_\ell,y_\ell),
\qquad \ell=1,\ldots,128\cdot512,
\]
enumerate the nodes of the computational grid.
The value $A_n(\boldsymbol{x}_\ell)$ is the stored complex
envelope at that node at the beginning of the slab.
We select nodes independently, with replacement, using
the probabilities
\[
p_\ell
=
\frac{|A_n(\boldsymbol{x}_\ell)|^2+10^{-16}}
{\displaystyle
	\sum_{m=1}^{128\cdot512}
	\bigl(|A_n(\boldsymbol{x}_m)|^2+10^{-16}\bigr)}.
\]

The computational grid spacings are chosen as follows
\[
\Delta x=\frac{96}{128}=0.75,
\qquad
\Delta y=\frac{304}{512}=0.59375.
\]
After selecting a node, we draw independent offsets
\[
\delta_x\sim\operatorname{Uniform}
\left[-\frac{\Delta x}{2},\frac{\Delta x}{2}\right],
\qquad
\delta_y\sim\operatorname{Uniform}
\left[-\frac{\Delta y}{2},\frac{\Delta y}{2}\right].
\]
The resulting sample coordinates are
\[
\begin{aligned}
	x&=\min\{48,\max\{-48,x_\ell+\delta_x\}\},\\
	y&=\min\{152,\max\{-152,y_\ell+\delta_y\}\}.
\end{aligned}
\]
 .

For complex sample values $q_1,\ldots,q_M$, we define
\[
\operatorname{mse}(q)
=\frac1M\sum_{r=1}^{M}|q_r|^2.
\]

The normalized neural envelope and physical wave function are now
\[
a_\theta(t,x,y)=\frac{A_\theta(t,x,y)}{0.06},
\qquad
\psi_\theta(t,x,y)
=0.06\,e^{\ii(-x+4y-6.9t)}a_\theta(t,x,y).
\]
For TDPSF, the normalized envelope residual is then
\[
\begin{aligned}
	R_\theta^{\rm TDPSF}
	={}&\ii\bigl(
	\partial_t a_\theta
	+0.6\,\partial_x a_\theta
	+3.6\,\partial_y a_\theta
	\bigr)\\
	&+\tfrac12\bigl(
	\partial_{xx}a_\theta
	+0.8\,\partial_{xy}a_\theta
	+\partial_{yy}a_\theta
	\bigr)
	+0.06^2|a_\theta|^2a_\theta.
\end{aligned}
\]
For PML, using the stretching operators $D_x$ and $D_y$
defined above, the corresponding residual is computed by the formula
\[
\begin{aligned}
	R_\theta^{\rm PML}
	={}&\frac{e^{-\ii(-x+4y-6.9t)}}{0.06}
	\Bigl[
	\ii\partial_t\psi_\theta+|\psi_\theta|^2\psi_\theta\\
	&\qquad
	+\tfrac12
	\bigl\{D_x^2+D_y^2+0.4(D_xD_y+D_yD_x)\bigr\}
	\psi_\theta
	\Bigr].
\end{aligned}
\]
 
We use $$R_\theta=R_\theta^{\rm TDPSF}$$ or
$$R_\theta=R_\theta^{\rm PML}$$ according to the boundary method.

If the $1536$ PDE samples are
$(t_n+\tau_r,x_r,y_r)$, the PDE contribution to the loss is
\[
\operatorname{mse}(R_\theta)
=
\frac1{1536}\sum_{r=1}^{1536}
\left|R_\theta(t_n+\tau_r,x_r,y_r)\right|^2.
\]

For the boundary contribution, we use $128$ matched
sample pairs on each pair of opposite faces.
For the left and right faces, independently sample
\[
t_r^x\sim\operatorname{Uniform}[t_n,t_{n+1}],
\qquad
y_r^B\sim\operatorname{Uniform}[-152,152],
\qquad r=1,\ldots,128.
\]
Each sample gives the two points
\[
(t_r^x,-48,y_r^B),
\qquad
(t_r^x,48,y_r^B).
\]
For the bottom and top faces, independently sample
\[
t_r^y\sim\operatorname{Uniform}[t_n,t_{n+1}],
\qquad
x_r^B\sim\operatorname{Uniform}[-48,48],
\]
giving the paired points
\[
(t_r^y,x_r^B,-152),
\qquad
(t_r^y,x_r^B,152).
\]
 
For TDPSF, the boundary loss enforces equality of the
envelope and its first coordinate derivative on opposite
faces. We have the formula:
\[
\begin{aligned}
	\mathcal L_B^{\rm TDPSF}
	=\frac1{128}\sum_{r=1}^{128}\Bigl[
	&|a_\theta(t_r^x,-48,y_r^B)
	-a_\theta(t_r^x,48,y_r^B)|^2\\
	+{}&
	|\partial_xa_\theta(t_r^x,-48,y_r^B)
	-\partial_xa_\theta(t_r^x,48,y_r^B)|^2\\
	+{}&
	|a_\theta(t_r^y,x_r^B,-152)
	-a_\theta(t_r^y,x_r^B,152)|^2\\
	+{}&
	|\partial_ya_\theta(t_r^y,x_r^B,-152)
	-\partial_ya_\theta(t_r^y,x_r^B,152)|^2
	\Bigr].
\end{aligned}
\]

In the case of PML, the boundary loss enforces zero envelope values
on the four outer faces:
$$
\begin{aligned}
	\mathcal L_B^{\rm PML}
	=\frac1{128}\sum_{r=1}^{128}\Bigl[
	&|a_\theta(t_r^x,-48,y_r^B)|^2
	+|a_\theta(t_r^x,48,y_r^B)|^2\\
	&+|a_\theta(t_r^y,x_r^B,-152)|^2
	+|a_\theta(t_r^y,x_r^B,152)|^2
	\Bigr].
\end{aligned}
$$
This is the sum of the four facewise mean squares.

The training objective is set as follows
\begin{equation}\label{eq:test1-loss}
	\mathcal L(\theta)
	=
	\operatorname{mse}(R_\theta)
	+5\mathcal L_B
	+10^{-7}\bigl(\|c_0\|_2^2+\|c_1\|_2^2\bigr),
\end{equation}
where $\mathcal L_B$ is the boundary loss for the selected
method.

The vectors $c_0,c_1\in\mathbb R^{384}$ collect, in a fixed
order, the real and imaginary parts of the complex output
coefficients $c_{\ell,a,b}^{(0)}$ and
$c_{\ell,a,b}^{(1)}$ defined in the network representation.
 
Consequently, we find
\[
\|c_0\|_2^2+\|c_1\|_2^2
=
\sum_{\ell=1}^{32}
\sum_{\substack{a,b\in\{0,1,2\}\\a+b\leq2}}
\left(
|c_{\ell,a,b}^{(0)}|^2
+|c_{\ell,a,b}^{(1)}|^2
\right).
\]
 
The transported initial term gives
$A_\theta(t_n,x,y)=A_n(x,y)$ exactly.

We describe the training on each slab  as follows.
First, 3 Gauss--Newton updates adjust $c_0$ and $c_1$
while the Gaussian parameters remain fixed.
Each update linearizes the sampled residuals with respect
to these output coefficients and solves the resulting
regularized least-squares problem.
Next, the limited-memory Broyden--Fletcher--Goldfarb--Shanno
method (L-BFGS) updates all $1024$ trainable parameters,
using at most $40$ iterations, or $120$ on the first slab.
Finally, two further Gauss--Newton updates adjust the
output coefficients.

The linear least-squares solves use a relative rank
tolerance of $10^{-5}$ to identify numerically independent
coefficient combinations. 

We use the three random seeds $0$, $1$, and $2$ for each
boundary method, giving 6 primary runs.
For each seed, the methods share the network architecture,
sampling rules and sample counts, loss weights, and
iteration caps.  

Accuracy is evaluated on the observation grid
$$
x_i=-24+i h_x,\qquad i=0,\ldots,64,
\qquad h_x=\frac{48}{64}=0.75,
$$
$$
y_j=-128+j h_y,\qquad j=0,\ldots,256,
\qquad h_y=\frac{256}{256}=1.
$$
As thus the observation grid contains $65\times257$ points.
Spatial $L^2$ norms are evaluated by the product
trapezoidal rule, as written below
\[
\|f\|_{L^2(\Omega_{\rm obs})}^2
\approx
h_xh_y
\sum_{i=0}^{64}\sum_{j=0}^{256}
w_i^x w_j^y |f(x_i,y_j)|^2,
\]
where
\[
w_i^x=
\begin{cases}
	\tfrac12,&i=0\text{ or }64,\\
	1,&1\leq i\leq63,
\end{cases}
\qquad \mbox{ and}
w_j^y=
\begin{cases}
	\tfrac12,&j=0\text{ or }256,\\
	1,&1\leq j\leq255.
\end{cases}
\]
The same quadrature evaluates the error numerator and
the fixed normalization
\[
S_{\rm ref}=\|\psi(0)\|_{L^2(\Omega_{\rm obs})}.
\]
Errors are recorded at $t_n=0.5n$, for $n=0,\ldots,200$.

The independent reference is computed in coordinates
moving with the carrier group velocity:
\[
\xi_x=x-4-0.6t,\qquad
\xi_y=y+72-3.6t.
\]
We define the envelope in these coordinates by
\[
A_{\rm mov}(t,\xi_x,\xi_y)
=
A(t,\xi_x+4+0.6t,\;\xi_y-72+3.6t).
\]
It satisfies
\[
\ii\partial_t A_{\rm mov}
=
-\tfrac12
\left(
\partial_{\xi_x}^2
+0.8\partial_{\xi_x}\partial_{\xi_y}
+\partial_{\xi_y}^2
\right)A_{\rm mov}
-|A_{\rm mov}|^2A_{\rm mov},
\]
with initial data
\[
A_{\rm mov}(0,\xi_x,\xi_y)
=
0.06\exp\left[-\frac{\xi_x^2+\xi_y^2}{72}\right].
\]

To evaluate the reference at an observation point
$(x_i,y_j)$, we first map that point to
$(\xi_x,\xi_y)$ and use cubic interpolation of the
moving-frame grid values.
When the mapped point lies outside the reference box,
the sampled envelope is assigned the value zero.
The carrier factor then reconstructs the physical
reference field $\psi_{\rm ref}$.
This independently computed field is used to assess
the PINN accuracy.

Reference sensitivity is assessed by comparing the
baseline reference with 3 alternative calculations.
Let $H_{\rm ref}$ denote the box half-width,
$N_{\rm ref}$ the number of grid points in each direction,
and $\delta t$ the reference time step.
The baseline settings are
\[
(H_{\rm ref},N_{\rm ref},\delta t)
=(128,512,0.01).
\]
For an alternative reference, define
\[
D(t)
=
\frac{
	\|\psi_{\rm ref,alt}(t)-\psi_{\rm ref,base}(t)\|
	_{L^2(\Omega_{\rm obs})}
}{S_{\rm ref}}.
\]
Here $\psi_{\rm ref,base}$ and $\psi_{\rm ref,alt}$
are the baseline and alternative reference fields,
evaluated at the same observation points and times.

The settings and measured differences are
\[
\begin{array}{l|ccc|cc}
	\text{variation}
	&H_{\rm ref}&N_{\rm ref}&\delta t&D(10)&D(50)\\ \hline
	\text{spatial grid}
	&128&768&0.01
	&2.02\times10^{-7}&1.15\times10^{-7}\\
	\text{box size}
	&192&768&0.01
	&9.74\times10^{-14}&2.60\times10^{-12}\\
	\text{time step}
	&128&512&0.005
	&6.05\times10^{-9}&2.00\times10^{-9}
\end{array}
\]

Finally, let
\[
\Bw=\Omega_{\rm comp}\setminus\Omega_{\rm obs}
\]
denote the buffer region.
For the computed PINN field $\psi_h$, the run records
contain the endpoint quantities
\[
\|\psi_h(t_n)\|_{L^2(\Omega_{\rm obs})},
\qquad
\|\psi_h(t_n)\|_{L^2(\Omega_{\rm comp})},
\qquad
\int_{\Bw}|\psi_h(t_n,x,y)|^4\,dx\,dy.
\]
These integrals are evaluated using computational-grid
sums with cell area $\Delta x\Delta y$.

\subsubsection{Numerical results for Test 1}
Figure~\ref{fig:test1-accuracy} gives the complete fixed-scale error
histories through $T=100$. At $t=10$, the median errors are
$2.98\times10^{-3}$ for TDPSF--PINN and $2.99\times10^{-3}$ for
PML--PINN. As the packet interacts with the layers, the TDPSF--PINN median
error is $2.28\times10^{-2}$ at $t=32$ and $3.85\times10^{-2}$ at
$t=50$ while the PML values are $1.29\times10^{-1}$ and
$8.35\times10^{-1}$ at $t=32$ and $50$. At $t=100$, we    can see that TDPSF reaches $4.35\times10^{-4}$,
compared with $2.75$ for PML.

This observation demonstrates the
advantage of  TDPSF--PINN for this opposed-velocity packet
at the shared Gaussian representation and training budget. In the long time simulation,  TDPSF--PINN is more accurate and stable.

In the picture, the lines represent  pointwise medians and bands span
the three paired seeds.
\space

\begin{figure}[htbp]\centering\color{black}\captionsetup{font={color=black}}
\includegraphics[width=\linewidth]{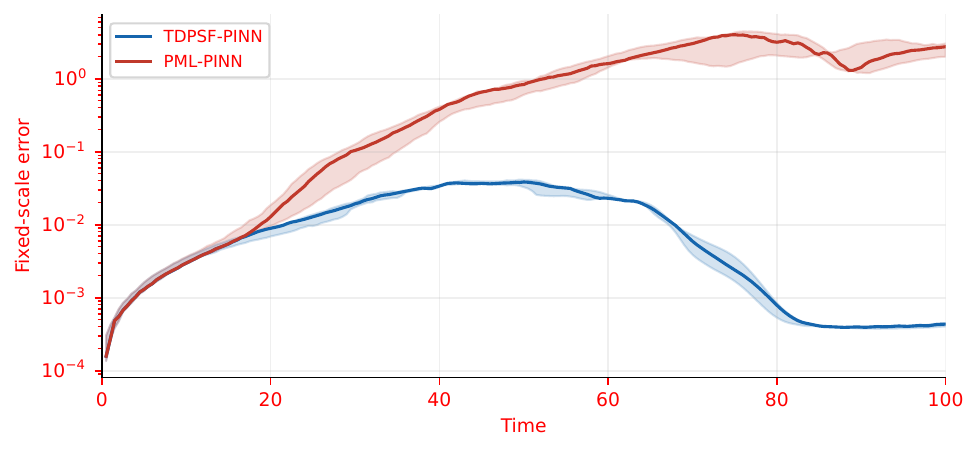}
\caption{Test 1: interior fixed-scale error versus time for TDPSF--PINN and PML--PINN through $T=100$. In the picture, the curves are pointwise medians and the bands give the full range over 3 paired seeds.  TDPSF--PINN is more stable and accurate in the long time run.}
\label{fig:test1-accuracy}
\end{figure}

\subsection{Test 2}
\label{sec:results-acoustic-crbc}
  We compare TDPSF--PINN and CRBC--PINN
through   $T_{\max}=32$.

\subsubsection{Acoustic model and boundary implementation}
  For $$U=(p,v_1,v_2)^T\colon[0,T_{\max}]\times\mathbb R^2
\longrightarrow\mathbb R^3,$$ we consider the equation  
\begin{equation}\label{eq:rev-acoustic}
 \partial_tU+A_1\partial_{x_1}U+A_2\partial_{x_2}U=0,
 \quad
 A_1=\begin{pmatrix}0&1&0\\1&0&0\\0&0&0\end{pmatrix},\qquad
 A_2=\begin{pmatrix}0&0&1\\0&0&0\\1&0&0\end{pmatrix}.
\end{equation}
Using the function $$\chi_{a,b}(x)=1-\vartheta((|x|-a)/(b-a)),$$  with the common quintic step $\vartheta$, we set the initial pressure  to be 
\begin{equation}\label{eq:rev-acoustic-initial}
\begin{aligned}
 p_0(x_1,x_2)
 &=C_0 e^{-((x_1+0.1)/0.7)^2-((x_2-0.05)/0.6)^2}
 \chi_{1,1.8}(x_1)\chi_{1,1.8}(x_2)\\
 &\qquad\times\cos\!\bigl(2.4(x_1+0.1)+3.2(x_2-0.05)\bigr),
 \qquad v_0=(0,0)^{\mathsf T}.
\end{aligned}
\end{equation}
The   constant $C_0$ is chosen such that

\[
\begin{aligned}
 \int_{\mathbb R^2}p_0(x)^2\,dx
 &=\int_{\mathbb R^2}
 e^{-2((x_1+0.1)/0.7)^2-2((x_2-0.05)/0.6)^2}
 \chi_{1,1.8}(x_1)^2\chi_{1,1.8}(x_2)^2\,dx\\
 &\approx0.6593215061578706.
\end{aligned}
\]
In this test, the initial data and this
normalization are identical for the 2  methods and the
independent reference calculation.
The independent whole-space reference has complementary Fourier
and Green-kernel representations. The
general Poisson representation of the pressure is defined as 
\begin{equation}\label{eq:test2-poisson-reference}
 p(t,x)=\frac{\partial}{\partial t}
 \left[
 \frac1{2\pi}\int_{|x-y|<t}
 \frac{p_0(y)}{\sqrt{t^2-|x-y|^2}}\,dy
 \right].
\end{equation}
Let 
$$
 \widehat p_0(k)=\int_{\mathbb R^2}p_0(y)e^{-\ii k\cdot y}\,dy,
 \qquad k=(k_1,k_2).
$$

The Fourier evolution factors are defined to be 
\[
 \widehat p(t,k)=\cos(|k|t)\widehat p_0(k),\qquad
 \widehat v(t,k)=-\ii k\frac{\sin(|k|t)}{|k|}\widehat p_0(k),
 \qquad k\neq0.
\]
At $k=0$, the pressure coefficient is $\widehat p_0(0)$ and
the velocity coefficient is zero. 

These factors are then evaluated on an
enlarged Fourier grid for $0\leq t<16$, before periodic images 
reach the evaluation region. 

Once the entire initial support lies strictly inside the backward
light cone, which  is,
\[
 t>\sup_{y\in\operatorname{supp}p_0}|x-y|,
\]
differentiation over the fixed initial support gives  
\begin{equation}\label{eq:rev-acoustic-reference}
\begin{aligned}
 p(t,x)&=-\frac{t}{2\pi}\int_{\operatorname{supp}p_0}
 \frac{p_0(y)}{(t^2-|x-y|^2)^{3/2}}\,dy,\\
 v(t,x)&=-\frac1{2\pi}\int_{\operatorname{supp}p_0}
 \frac{(x-y)p_0(y)}{(t^2-|x-y|^2)^{3/2}}\,dy.
\end{aligned}
\end{equation}

The early Fourier calculation has $2048^2$ points on
$[-16,16)^2$. For $t\geq16$, we split each source coordinate into
$[-1.8,-1]$, $[-1,1]$, and $[1,1.8]$, apply $48$ Gaussian quadrature
nodes on each panel, and retain moment orders $$0,\ldots,24$$ in our
Green-kernel expansion. The strict inequalities defining
$$\Omega_{\rm obs}=(-1.5,1.5)^2$$ and
$$\Omega_{\rm phys}=(-2,2)^2$$ select $47^2$ and $63^2$ nodes,
respectively, from the $128^2$ endpoint grid on $[-4,4)^2$.

\begin{remark}
    In Test 2, choosing $\Omega_{\rm obs}$   entirely inside $\Omega_{\rm phys}$ gives a focused accuracy measurement in the central region, separated from the buffer interface by a distance of \(0.5\). This choice is useful when  the application cares specifically about the solution near the center. 
\end{remark}

We check the reference independently of the neural training.
At $t=1$, increasing the Fourier grid from $2048^2$ to $4096^2$
gives a fixed-scale $L^2(\Omega_{\rm obs})$ difference of
$3.44\times10^{-9}$. At $t=16$, using $96$ Gaussian nodes per
panel and moment orders $$0,\ldots,36$$ gives a largest componentwise
difference of $6.02\times10^{-18}$ on the endpoint grid.
The accompanying reference-validation code records the full-grid
comparisons at $$t=0.25,1,4,7,12,15.75$$ and the late-time comparisons
at $t=16,32$.

\paragraph{\bf CRBC equations.}
We use the complete radiation boundary condition (CRBC)
construction of Hagstrom and Warburton in 
\cite{HagstromWarburton2009}.
The following definitions specify the auxiliary variables,
their signs, and their coupling at the corners.

The computational square, physical square, and observation
square are set to be
\[
\Omega_{\rm comp}=(-H,H)^2,\qquad H=4,
\]

and 
$$
\Omega_{\rm phys}=(-2,2)^2,
\qquad
\Omega_{\rm obs}=(-1.5,1.5)^2.
$$

We number the four outer edges by $$e=1,\ldots,4.$$

For each edge, $\boldsymbol{x}_e(z)$ denotes its spatial
parameterization, where $z\in[-H,H]$ increases in the
positive Cartesian tangential direction.
The outward unit normal is denoted by $n_e$.
The normal and tangential matrices are specified as follows

\[
\begin{array}{c|c|c|c|c|c}
	e&\text{edge}&\boldsymbol{x}_e(z)&n_e
	&A_{n_e}&A_{{\rm tan},e}\\ \hline
	1&\text{right}&(H,z)&(1,0)&A_1&A_2\\
	2&\text{left}&(-H,z)&(-1,0)&-A_1&A_2\\
	3&\text{top}&(z,H)&(0,1)&A_2&A_1\\
	4&\text{bottom}&(z,-H)&(0,-1)&-A_2&A_1
\end{array}
\]

Thus, we can write
\[
A_{n_e}=n_{e,1}A_1+n_{e,2}A_2.
\]

On edge $e$, we introduce the 3-component fields
\[
W_{e,j}(t,z)\in\mathbb R^3,
\qquad j=0,\ldots,P.
\]

At level zero, we can see that
\[
W_{e,0}(t,z)=U(t,\boldsymbol{x}_e(z)),
\]
which is the physical field evaluated on the edge.

For a fixed edge, we suppress $e$ and write
\[
W_j=(p_j,v_{1,j},v_{2,j})^{\mathsf T}.
\]
The normal velocity component is then
\[
v_{n,j}=n_1v_{1,j}+n_2v_{2,j},
\]
and the tangential velocity component is as follows
\[
v_{{\rm tan},j}
=
\begin{cases}
	v_{2,j},&\text{on the left and right edges},\\
	v_{1,j},&\text{on the bottom and top edges}.
\end{cases}
\]

For $P>1$, the prescribed scalar parameters are now
\begin{equation}\label{eq:rev-crbc-parameters}
	a_j=a_{\min}^{j/(P-1)},
	\qquad
	b_j=\frac{1-a_j^2}{a_jT_{\max}},
	\qquad j=0,\ldots,P-1.
\end{equation}
The time operator $S_j$ acts on a scalar or vector field
$f$ by
\[
S_jf=a_j\partial_t f+b_jf,
\]
with componentwise action for a vector field.
For $P=1$, we use $a_0=1$ and $b_0=0$.

In the reported calculations,
$$
P=4,\qquad a_{\min}=0.05,\qquad T_{\max}=32.
$$

For each edge, we define the vector residual
\begin{equation}\label{eq:rev-crbc-edge}
	\begin{aligned}
		(R_e)_j
		={}&
		\partial_t(W_{e,j}+W_{e,j+1})\\
		&+A_{{\rm tan},e}\partial_z(W_{e,j}+W_{e,j+1})\\
		&+A_{n_e}S_j(W_{e,j+1}-W_{e,j}),
		\qquad j=0,\ldots,P-1.
	\end{aligned}
\end{equation}
The edge equations require $(R_e)_j=0$.

At every level, including the physical level $j=0$,
we also impose the individual tangential equation
\[
(T_e)_j
:=
\partial_t v_{{\rm tan},e,j}
+\partial_z p_{e,j}
=0,
\qquad j=0,\ldots,P,
\]
in which $p_{e,j}$ and $v_{{\rm tan},e,j}$ are the pressure
and tangential-velocity components of $W_{e,j}$.

For a vector $W=(W_p,W_{v_1},W_{v_2})^{\mathsf T}$,
define
\[
C_{n_e}^-W
=
W_p-n_{e,1}W_{v_1}-n_{e,2}W_{v_2}.
\]
 
The last auxiliary level satisfies
\[
C_e:=C_{n_e}^-W_{e,P}=0.
\]

At a corner, let
\[
c=(s_1,s_2)\in
\mathcal C
:=
\{(1,1),(1,-1),(-1,1),(-1,-1)\}.
\]
The corner location is $$(s_1H,s_2H).$$

At corner $c$, we introduce
\[
W_{j,k}^{(c)}(t)\in\mathbb R^3,
\qquad 0\leq j,k\leq P.
\]

The values $W_{j,0}^{(c)}$ are supplied by the vertical
edge fields evaluated at tangential coordinate $z=s_2H$.
The values $W_{0,k}^{(c)}$ are supplied by the horizontal
edge fields at $z=s_1H$.
In particular,
$$
W_{0,0}^{(c)}(t)=U(t,s_1H,s_2H).
$$

For a fixed corner, we suppress the superscript $(c)$
and write
\[
W_{j,k}=(p_{j,k},u_{j,k},v_{j,k})^{\mathsf T},
\]
where $u_{j,k}$ and $v_{j,k}$ are the auxiliary velocity
components in the $x_1$ and $x_2$ directions.

For $0\leq j,k<P$, we define the pressure residual
\begin{align}\label{eq:rev-crbc-corner}
	(R_c^p)_{j,k}
	={}&
	\partial_t
	(p_{j+1,k+1}+p_{j,k+1}+p_{j+1,k}+p_{j,k})
	\nonumber\\
	&+s_1S_j
	(u_{j+1,k+1}-u_{j,k+1}+u_{j+1,k}-u_{j,k})
	\nonumber\\
	&+s_2S_k
	(v_{j+1,k+1}-v_{j+1,k}+v_{j,k+1}-v_{j,k}).
\end{align}
The pressure equation is then $$(R_c^p)_{j,k}=0.$$

The velocity residuals are then
\begin{align}\label{eq:rev-crbc-corner-velocity}
	(R_c^u)_{j,k}
	&=
	\partial_t(u_{j+1,k}+u_{j,k})
	+s_1S_j(p_{j+1,k}-p_{j,k}),
	\nonumber\\[-1mm]
	&\hspace{35mm}0\leq j<P,\quad 0\leq k\leq P,
	\\[1mm]
	(R_c^v)_{j,k}
	&=
	\partial_t(v_{j,k+1}+v_{j,k})
	+s_2S_k(p_{j,k+1}-p_{j,k}),
	\nonumber\\[-1mm]
	&\hspace{35mm}0\leq j\leq P,\quad 0\leq k<P.
	\nonumber
\end{align}
We impose both $(R_c^u)_{j,k}=0$ and
$$(R_c^v)_{j,k}=0$$ over their stated index ranges.

The terminal residuals are
\[
(C_c^x)_k=p_{P,k}-s_1u_{P,k},
\qquad k=0,\ldots,P-1,
\]
\[
(C_c^y)_j=p_{j,P}-s_2v_{j,P},
\qquad j=0,\ldots,P-1,
\]
and
\begin{equation}\label{eq:rev-crbc-terminal-corner}
	C_c^{xy}
	=
	2p_{P,P}-s_1u_{P,P}-s_2v_{P,P}.
\end{equation}
All three terminal residuals are set to zero.
 
These edge and corner equations specify the component
CRBC system used in the PINN, following
\cite[Eqs.~(7.7)--(7.14)]{HagstromWarburton2009}.
At $t=0$, all edge auxiliary fields with $j\geq1$ and
all mixed corner fields with $j,k\geq1$ are zero.

The physical boundary trace also initially vanishes 
as the initial pulse is supported away from the
outer boundary.
At later times, the computed auxiliary states provide
the initial data for the next slab of computations.

\paragraph{\bf TDPSF filtering and diagnostics.}
For a nonzero Fourier wave vector $k=(k_1,k_2)$, define
\[
|k|=\sqrt{k_1^2+k_2^2},
\qquad
d(k)=\frac{k}{|k|}.
\]

The 2 propagating acoustic branches, indexed by
$\sigma\in\{-1,1\}$, have frequencies and group velocities
\[
\omega_\sigma(k)=\sigma|k|,
\qquad
v_{g,\sigma}(k)=\nabla_k\omega_\sigma(k)=\sigma d(k).
\]

For a Fourier coefficient
$$\widehat U=(\widehat p,\widehat v_1,\widehat v_2)^{\mathsf T},$$ 
we write $\widehat v=(\widehat v_1,\widehat v_2)$.
The projector onto branch $\sigma$ is as thus
\[
\Pi_\sigma(k)\widehat U
=
\frac12
\bigl(\widehat p+\sigma d(k)\cdot\widehat v\bigr)
\begin{pmatrix}
	1\\
	\sigma d_1(k)\\
	\sigma d_2(k)
\end{pmatrix}.
\]
 
The buffer width is now
\[
w=H-2=2.
\]
For edge $e$, we define the outward distance from the
physical square by
\[
\zeta_e(x)=n_e\cdot x-2.
\]
 
The four distances, in application order, are
\[
\zeta_1=x_1-2,\qquad
\zeta_2=-x_1-2,\qquad
\zeta_3=x_2-2,\qquad
\zeta_4=-x_2-2.
\]
The spatial windows are
\[
\eta_e(x)
=
\vartheta\!\left(\frac{\zeta_e(x)}{0.9w}\right)
\vartheta\!\left(\frac{w-\zeta_e(x)}{0.1w}\right),
\qquad e=1,\ldots,4,
\]
where $\vartheta$ is the quintic step defined above.

For branch $\sigma$, the outward normal group velocity
at edge $e$ is $n_e\cdot v_{g,\sigma}(k)$.
Its removal weight is
\[
m_{e,\sigma}(k)
=
\left[
1+\exp\left(
-\frac{n_e\cdot v_{g,\sigma}(k)-0.2}{0.1}
\right)
\right]^{-1},
\qquad k\neq0,
\]
with $m_{e,\sigma}(0)=0$.

Let $U^{\rm raw}$ be the endpoint values of the field
network on the computational grid before filtering.
Starting from $U^{[0]}=U^{\rm raw}$, we apply
\[
U^{[e]}
=
U^{[e-1]}
-
\eta_e\mathcal F^{-1}
\left[
\sum_{\sigma\in\{-1,1\}}
m_{e,\sigma}\Pi_\sigma
\mathcal F\bigl(\eta_e U^{[e-1]}\bigr)
\right],
\qquad e=1,\ldots,4.
\]
in which $\mathcal F$ and $\mathcal F^{-1}$ are the
componentwise discrete Fourier transform and its inverse.
 
The final real-valued endpoint is
\[
U^{\rm filt}=\operatorname{Re}U^{[4]}.
\]

For completeness, the scalar spatial cutoff used in
the initial data is defined as follows
\[
\chi_{a,b}(s)
=
1-\vartheta\!\left(\frac{|s|-a}{b-a}\right),
\qquad 0\leq a<b.
\]
 
For any spatial region $D$, the acoustic energy norm is
\begin{equation}\label{eq:rev-acoustic-norm}
	\|U(t)\|_{E(D)}^2
	=
	\int_D
	\bigl(p(t,x)^2+v_1(t,x)^2+v_2(t,x)^2\bigr)\,dx.
\end{equation}
With unit density and sound speed, the physical energy
in $D$ is $\tfrac12\|U(t)\|_{E(D)}^2$.
We therefore record
\[
\mathcal E_{\rm phys}(t)
=
\tfrac12\|U(t)\|_{E(\Omega_{\rm phys})}^2,
\qquad
\mathcal E_{\rm buffer}(t)
=
\tfrac12\|U(t)\|_{E(\mathcal B_w)}^2,
\]
where
$
\mathcal B_w
=
\Omega_{\rm comp}\setminus\Omega_{\rm phys}.
$

The CRBC auxiliary diagnostics are then
\[
\mathcal A_{\rm edge}(t)^2
=
\sum_{e=1}^{4}\sum_{j=1}^{P}
\int_{-H}^{H}|W_{e,j}(t,z)|_{\mathbb R^3}^2\,dz,
\]
and
\[
\mathcal A_{\rm corner}(t)^2
=
\sum_{c\in\mathcal C}
\sum_{j=1}^{P}\sum_{k=1}^{P}
|W_{j,k}^{(c)}(t)|_{\mathbb R^3}^2,
\]
where

$$|(q_1,q_2,q_3)|_{\mathbb R^3}^2=q_1^2+q_2^2+q_3^2.$$
 
The prescribed initial-value problem is propagated
through $t=8$, $16$, and $32$.
These are the times $T_0$, $2T_0$, and $4T_0$, with $T_0=8$
and final horizon $T_{\max}=32$.

\paragraph{\bf Neural and filter parameters.}
Let
\[
t_n=n\Delta t,\qquad \Delta t=0.25,
\qquad \tau=t-t_n\in[0,\Delta t].
\]
The physical field network is denoted by
\[
U_\theta=(p_\theta,v_{1,\theta},v_{2,\theta})^{\mathsf T},
\]
where $\theta$ collects its trainable weights and biases.

Its four input coordinates are
$$
q(t,x)
=
\begin{pmatrix}
	\tau\\
	x_1/4\\
	x_2/4\\
	(r(x)-t)/4
\end{pmatrix},
\qquad
r(x)=
\sqrt{(x_1+0.1)^2+(x_2-0.05)^2+0.01}.
$$

The field network has two hidden layers, each containing
$48$ neurons, followed by a linear layer with 3 outputs.

Both hidden layers has the scalar activation function
\[
\sigma_{\rm act}(s)=\sin(15s),
\]
where $s$ is the  input to a neuron.
Including biases, the field parameter count is then 
$2739$.

In the case of CRBC, each edge has its own network with inputs
$(\tau,z/H)$ and outputs
$$
(W_{e,1},\ldots,W_{e,P}),
$$
interpreted as a vector of $12$ scalar components.
Each corner has its own network with input $\tau$ and
outputs the $48$ components of
$W_{j,k}^{(c)}$, $j,k=1,\ldots,P$.

Each auxiliary network has 2 hidden layers of width
$24$ and uses the same sinusoidal activation.
An edge network contains $972
$
parameters, and a corner network contains
$
1848
$
parameters.

The 4 edge and four corner networks therefore contain
$11280
$
auxiliary parameters.

On each slab, the primary calculation uses
$
N_P=1024,\qquad N_I=1024,\qquad N_B=16.
$
Here $N_P$ is the number of space--time PDE samples,
$N_I$ is the number of physical-field initial samples,
and $N_B$ is the number of boundary samples per edge
or corner.

The PDE samples are uniform in
$$[t_n,t_{n+1}]\times\Omega_{\rm comp}.$$
The initial samples lie at $t=t_n$:
$512$ spatial points are uniform in $\Omega_{\rm comp}$
and 512 are uniform in $\Omega_{\rm phys}$.

For boundary residuals, we draw $N_B$ pairs
$(t_r,z_r)$ uniformly from
$$[t_n,t_{n+1}]\times[-H,H]$$ and evaluate the four
edge parameterizations $\boldsymbol{x}_e(z_r)$.
The corner residuals use the same $N_B$ sample times.
Initial matching of the edge auxiliaries uses $16$
spatial samples on each edge at $t_n$.
Each corner's initial matching uses its stored
$$P\times P\times3$$ state at $t_n$.

The endpoint field is evaluated on a $128\times128$
grid, with spacing
$$
\Delta x_1=\Delta x_2= \frac1{16}.
$$
 
For CRBC, cubic interpolation also represents the
preceding edge states. The corner states are stored
directly as vectors.
In the case of TDPSF, the initial interpolation uses the filtered
endpoint.

In our first slab, 10 Adam steps precede at most
$900$ L-BFGS iterations.
Each subsequent slab then uses at most $300$ L-BFGS iterations.
All field and auxiliary weights are trained in those slabs.

In our implementation, when advancing to a new slab, the first network input
changes from the old local time to the new local time.
For any of the networks, let $B^{(1)}$ and $b^{(1)}$
be its first-layer weight matrix and bias.
The carried parameters are shifted by
$$
b^{(1)}_{\rm new}
=
b^{(1)}_{\rm old}
+\Delta t\,B^{(1)}_{\rm old}(:,1),
$$
where $B^{(1)}(:,1)$ is the column multiplying local time.

We also note that other carried weights are unchanged by this
coordinate adjustment.

Filtering is performed after every TDPSF slab, with 
$$\Delta t_{\rm filter}=0.25.$$

The primary comparison uses the paired seeds
$0,1,2,3,4$.
The joint refinement uses seed $0$ for both methods,
with the following settings:
$$
\begin{array}{l|cc}
	&\text{primary}&\text{joint refinement}\\ \hline
	\text{field hidden width}&48&64\\
	\text{auxiliary hidden width}&24&32\\
	N_P&1024&2048\\
	N_I&1024&2048\\
	N_B&16&16\\
	\text{ordinary-slab L-BFGS cap}&300&600\\
	\text{first-slab L-BFGS cap}&900&1800
\end{array}
$$

\paragraph{\bf Exact loss normalization and endpoint handover.}
For any sampled array $Q$ containing $N_Q$ scalar entries, we 
define
\[
\operatorname{mse}(Q)
=
\frac1{N_Q}\sum_{\nu=1}^{N_Q}|Q_\nu|^2.
\]

The enumeration $\nu$ includes every sample, auxiliary
index, and component belonging to that array.

The 3-component acoustic PDE residual is
\[
\mathcal R_{\rm PDE}
=
\partial_tU_\theta
+A_1\partial_{x_1}U_\theta
+A_2\partial_{x_2}U_\theta
=
\begin{pmatrix}
	\partial_tp_\theta
	+\partial_{x_1}v_{1,\theta}
	+\partial_{x_2}v_{2,\theta}\\
	\partial_tv_{1,\theta}+\partial_{x_1}p_\theta\\
	\partial_tv_{2,\theta}+\partial_{x_2}p_\theta
\end{pmatrix}.
\]
Thus $\mathcal R_{\rm PDE}$ has $N_P\times3$ sampled
scalar entries.
The objective is then
\begin{equation}\label{eq:test2-loss}
	\mathcal L
	=
	\mathcal L_P+100\mathcal L_I+5\mathcal L_B,
	\qquad
	\mathcal L_P=3\operatorname{mse}(\mathcal R_{\rm PDE}).
\end{equation}
Equivalently, we have
\[
\mathcal L_P
=
\frac1{N_P}\sum_{r=1}^{N_P}
|\mathcal R_{\rm PDE}(t_r,x_r)|_{\mathbb R^3}^2.
\]

Let $U_n(x)$ denote the known interpolant of the stored
field endpoint at time $t_n$.
At an initial sample $x_r^I$, we  define
$$
(\delta U)_r
=
U_\theta(t_n,x_r^I)-U_n(x_r^I),
\qquad r=1,\ldots,N_I.
$$
In the case of TDPSF, we have
\[
\mathcal L_I=3\operatorname{mse}(\delta U).
\]

For CRBC, let $W_{e,j}^{\,n}(z)$ be the known interpolant
of the preceding edge state and let
$W_{j,k}^{(c),n}$ be the stored corner state.
Using a subscript $\theta$ for the current neural
approximations, we define
\[
(\delta W_e)_{r,j}
=
W_{e,j,\theta}(t_n,z_r^I)
-W_{e,j}^{\,n}(z_r^I),
\qquad j=1,\ldots,P,
\]
and
$$
(\delta W_c)_{j,k}
=
W_{j,k,\theta}^{(c)}(t_n)
-W_{j,k}^{(c),n},
\qquad j,k=1,\ldots,P.
$$

The CRBC initial loss is
\begin{equation}\label{eq:test2-initial-loss}
	\mathcal L_I
	=
	3\operatorname{mse}(\delta U)
	+\frac18\sum_{e=1}^{4}\operatorname{mse}(\delta W_e)
	+\frac18\sum_{c\in\mathcal C}\operatorname{mse}(\delta W_c).
\end{equation}

For the TDPSF boundary loss, we define the paired traces
\[
\begin{aligned}
	U_{1,+}(t,z)&=U_\theta(t,H,z),&
	U_{1,-}(t,z)&=U_\theta(t,-H,z),\\
	U_{2,+}(t,z)&=U_\theta(t,z,H),&
	U_{2,-}(t,z)&=U_\theta(t,z,-H).
\end{aligned}
\]
 
The boundary loss is defined to be
\[
\mathcal L_B^{\rm TDPSF}
=
\frac32\left[
\operatorname{mse}(U_{1,+}-U_{1,-})
+
\operatorname{mse}(U_{2,+}-U_{2,-})
\right].
\]

In the case of CRBC, the symbols $R_e,T_e,C_e$ denote the sampled
edge residuals defined above.
The symbols $R_c^p,R_c^u,R_c^v$ and
$C_c^x,C_c^y,C_c^{xy}$ denote the sampled corner residuals.
Their array dimensions are given below
\[
\begin{array}{c|c}
	\text{residual array}&\text{dimensions}\\ \hline
	R_e&N_B\times P\times3\\
	T_e&N_B\times(P+1)\\
	C_e&N_B\\
	R_c^p&N_B\times P\times P\\
	R_c^u&N_B\times P\times(P+1)\\
	R_c^v&N_B\times(P+1)\times P\\
	C_c^x,\ C_c^y&N_B\times P\\
	C_c^{xy}&N_B
\end{array}
\]

For example, we have
\[
\operatorname{mse}(R_e)
=
\frac1{3N_BP}
\sum_{r=1}^{N_B}\sum_{j=0}^{P-1}\sum_{\mu=1}^{3}
|(R_e)_{r,j,\mu}|^2,
\]
where $\mu$ selects a vector component, whereas
\[
\operatorname{mse}(R_c^p)
=
\frac1{N_BP^2}
\sum_{r=1}^{N_B}\sum_{j=0}^{P-1}\sum_{k=0}^{P-1}
|(R_c^p)_{r,j,k}|^2.
\]

The implemented CRBC boundary loss is now
\begin{align}\label{eq:test2-boundary-loss}
	\mathcal L_B^{\rm CRBC}
	={}&
	\frac18\sum_{e=1}^{4}
	\left\{
	\operatorname{mse}(R_e)
	+\operatorname{mse}(T_e)
	+\operatorname{mse}(C_e)
	\right\}
	\nonumber\\
	&+\frac18\sum_{c\in\mathcal C}
	\left\{
	\frac{
		\operatorname{mse}(R_c^p)
		+\operatorname{mse}(R_c^u)
		+\operatorname{mse}(R_c^v)}{3}
	\right.
	\nonumber\\
	&\hspace{35mm}\left.
	+\operatorname{mse}(C_c^x)
	+\operatorname{mse}(C_c^y)
	+\operatorname{mse}(C_c^{xy})
	\right\}.
\end{align}

At the end of a slab, we let
\[
U_{n+1}^{\rm raw}(x)=U_\theta(t_{n+1},x)
\]
denote the neural endpoint evaluated on the grid.
CRBC stores this raw endpoint together with
\[
W_{e,j,\theta}(t_{n+1},z),
\qquad
W_{j,k,\theta}^{(c)}(t_{n+1}),
\]
for the auxiliary index ranges defined above.
TDPSF applies the four-window filter and stores
$U_{n+1}^{\rm filt}$.
These stored states then define the initial targets for
the next slab.

For saved-weight verification, we distinguish the physical endpoint time
from the time origin of the slab on which the weights were trained.
The final weights are evaluated at $t=32$ with local time origin
$31.75$. We then apply the TDPSF endpoint filter for TDPSF--PINN,
or retain the raw endpoint for CRBC--PINN. We reconstruct all twelve
final endpoints and compare the field components with their saved
values using relative tolerance $10^{-12}$ and absolute tolerance
$10^{-14}$.

\subsubsection{Numerical results for Test 2}

Figure~\ref{fig:test2-accuracy} gives the complete error histories through
$T_{\max}=32$.  
In our figure, the left panel gives pointwise medians and the full range over five
paired seeds, while the right panel gives the individual refined seed-0
runs.

  At $t=8$, the primary
median errors are $4.52\times10^{-2}$ and $9.33\times10^{-2}$;
at $t=32$, they are $1.01\times10^{-2}$ and $2.89\times10^{-1}$,
respectively for   TDPSF--PINN and
  CRBC--PINN. Our reported advantage of TDPSF--PINN concerns the
smaller residual errors and their sustained decay at late times.

We compare the joint enlargement with the primary run at the same
seed. For seed~0, the final TDPSF--PINN error changes from
$2.88\times10^{-2}$ to $7.16\times10^{-3}$, and the final
CRBC--PINN error changes from $2.21\times10^{-1}$ to
$1.00\times10^{-1}$.

The enlargement changes the network widths,
PDE and initial sample counts, and optimizer budgets together while
the boundary sample count, CRBC order, slab length, and endpoint grid
remain fixed.  

To assess the neural approximation between training points, we
evaluate the stored primary seed-0 first-slab networks and their PDE
residuals at $16384$ independent uniform points in the domain
$$[0,0.25]\times(-4,4)^2,$$ using NumPy's independent generator with
seed $62931$.

We average the squared Euclidean norm of the
3-component residual, using the same component normalization as
$\mathcal L_P$ in Equation~\eqref{eq:test2-loss}.

In this test, every primary TDPSF--PINN final error is below every primary
CRBC--PINN final error. Moreover, the ratio of the two primary medians is
approximately $28.5$. All 5 paired primary runs and both refined
runs complete the given horizon.  

\begin{figure}[!htbp]
\centering\color{black}\captionsetup{font={color=black}}
\includegraphics[width=\linewidth]{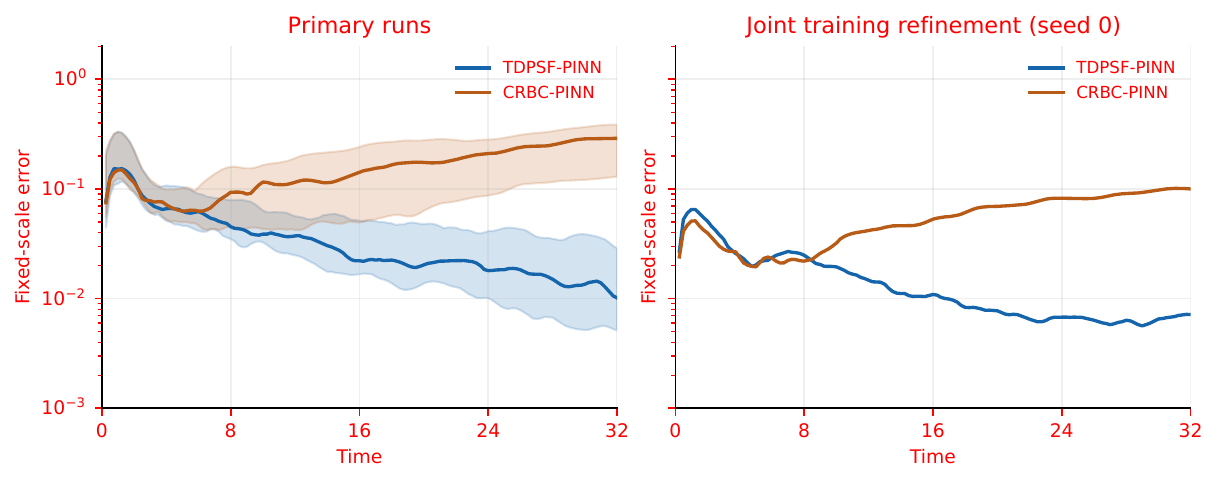}
\caption{Test 2: fixed-scale error versus time for the primary PINN runs (left) and joint training refinement at seed~0 (right). In the left panel, the curves are pointwise medians and the bands show the full range over five paired seeds. The right panel shows the individual refined runs with seed~0. Both panels cover the long time range $0\leq t\leq T_{\max}=32$.  }
\label{fig:test2-accuracy}
\end{figure}

\begin{figure}[H]
\centering
\includegraphics[width=\linewidth]{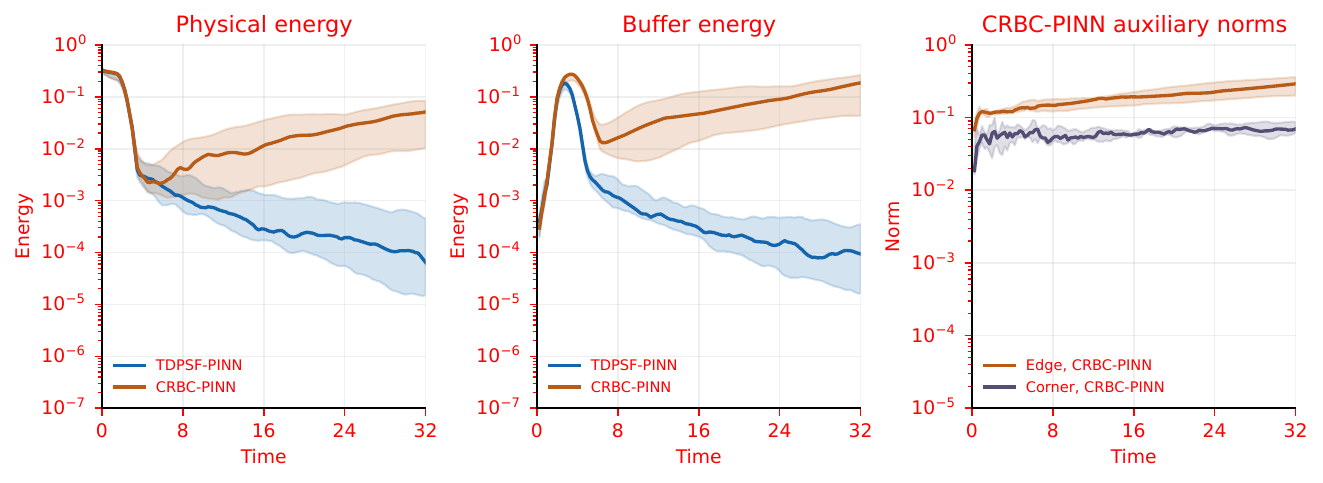}
\captionsetup{font={color=black}}
\caption{Test 2: physical energy, buffer energy, and CRBC edge/corner auxiliary norms. The energy convention is $\mathcal E_D=\tfrac12\int_D(p^2+|v|^2)\,dx$. The curves are pointwise medians and the bands show the full range over 5 seeds. The energy panels start at $10^{-7}$ and the auxiliary panel at $10^{-5}$.}
\label{fig:rev-acoustic-energy-auxiliary}
\end{figure}

Figure~\ref{fig:rev-acoustic-energy-auxiliary} shows the long-time
decay of the TDPSF--PINN energy in both the physical region and buffer.

 We can see that after the outgoing pulse passes through these regions, their energies
decrease overall, with small local increases. At $t=32$, the primary
median physical energies are $6.46\times10^{-5}$ for TDPSF--PINN
and $5.12\times10^{-2}$ for CRBC--PINN while the corresponding buffer
energies are $9.39\times10^{-5}$ and $1.85\times10^{-1}$.

Together with the smaller late-time errors in
Figure~\ref{fig:test2-accuracy},  those histories demonstrate the
advantage of TDPSF--PINN in eliminating late-time residual numerical
activity. We note that the primary runs use the shared field architecture and
ordinary-slab budget of $300$ L-BFGS iterations specified above.

The complete refined seed-0 trajectories use $1834.119$ process
CPU seconds for TDPSF--PINN and $5642.886$ for CRBC--PINN,
including training and endpoint operations.

\subsection{Test 3}
\label{sec:results-test3}
  
The comparison is between the time-dependent phase-space filter
combined with a physics-informed neural network (TDPSF--PINN) and
the second-order Bayliss--Turkel boundary condition combined with
a physics-informed neural network (BT2--PINN). Both computations
cover the time interval $0\leq t\leq100$.

\subsubsection{Acoustic model and boundary implementation}

 Let $x=(x_1,x_2)\in\mathbb R^2$ denote position and let $|x|$ denote
its Euclidean norm.

Obstacle $j$ ($j=1,2,3$)
is the closed disk
\[
D_j=\{x\in\mathbb R^2:|x-c_j|\leq a\},
\]
where $a$ is the common radius and $c_j$ is its center:
\begin{equation}\label{eq:test3-geometry}
	a=0.75,\qquad
	c_1=(-1,0),\qquad
	c_2=(1,0),\qquad
	c_3=(0.2,1.6),\qquad
	\Omega=\mathbb R^2\setminus\bigcup_{j=1}^3D_j.
\end{equation}
 
On the boundary $\partial D_j$, define
\[
\nu_j(x)=\frac{x-c_j}{a}.
\]

Let $p(t,x)\in\mathbb R$ denote acoustic pressure and let
$v(t,x)=(v_1(t,x),v_2(t,x))^T\in\mathbb R^2$ denote acoustic
velocity. We have the system
\begin{equation}\label{eq:test3-acoustics}
	\begin{aligned}
		\partial_t p+\nabla\cdot v&=q(x)g(t),
		&&x\in\Omega,\quad t>0,\\
		\partial_t v+\nabla p&=0,
		&&x\in\Omega,\quad t>0,\\
		p(0,x)&=0,\qquad v(0,x)=0,
		&&x\in\Omega,\\
		v(t,x)\cdot\nu_j(x)&=0,
		&&x\in\partial D_j,\quad j=1,2,3,
	\end{aligned}
\end{equation}
in which
\[
\nabla p=(\partial_{x_1}p,\partial_{x_2}p)^T,
\qquad
\nabla\cdot v=\partial_{x_1}v_1+\partial_{x_2}v_2.
\]
 
Equivalently, the pressure satisfies
\[
\partial_t^2p-\Delta p=q(x)g'(t),
\qquad
\partial_{\nu_j}p=0\quad\text{on }\partial D_j,
\qquad
p(0,x)=\partial_t p(0,x)=0,
\]
where
 $\partial_{\nu_j}p=\nu_j\cdot\nabla p$.

 Let $x_s$ denote the source center
and let $r_s$ denote its support radius. They are defined by the following formulas
\begin{align}
	q(x)&=
	\begin{cases}
		20\exp\!\left(
		1-\dfrac{1}{1-|x-x_s|^2/r_s^2}
		\right),&|x-x_s|<r_s,\\
		0,&|x-x_s|\geq r_s,
	\end{cases}
	\qquad
	x_s=(0,0.55),\qquad r_s=0.15,
	\label{eq:test3-source}\\
	g(t)&=
	\begin{cases}
		\sin^4(\pi t/8)
		\bigl\{\sin(2t)+\tfrac12\sin(4.5t)\bigr\},&0<t<8,\\
		0,&\text{otherwise}.
	\end{cases}
	\label{eq:test3-pulse}
\end{align}

We collect the pressure and velocity components into
$$
U=(p,v_1,v_2)^T.
$$
Both methods use the decomposition
\[
U=U^{\rm inc}+U^s,
\qquad
U^{\rm inc}=(p^{\rm inc},v_1^{\rm inc},v_2^{\rm inc})^T,
\qquad
U^s=(p^s,v_1^s,v_2^s)^T.
\]
The incident field $U^{\rm inc}$ is the solution generated by
the prescribed source in the whole plane. The scattered field
$U^s$ accounts for the obstacle interactions.

The incident calculation uses the distance
$$
\rho=|x-x_s|
$$
from the source center. 

We write 
\[
q_0(\rho)=q(x_s+(\rho,0)),
\qquad
\widehat q(\kappa)
=2\pi\int_0^{r_s}q_0(\rho)J_0(\kappa\rho)\rho\,d\rho,
\]
where $\kappa\geq0$ is the radial spatial wavenumber, and
$J_0$ and $J_1$ denote Bessel functions of the first kind
of orders zero and one. Define the temporal convolutions
\[
\begin{aligned}
	G_c(\kappa,t)
	&=\int_0^{\min(t,8)}
	g(s)\cos\!\bigl(\kappa(t-s)\bigr)\,ds,\\
	G_s(\kappa,t)
	&=\int_0^{\min(t,8)}
	g(s)\sin\!\bigl(\kappa(t-s)\bigr)\,ds.
\end{aligned}
\]
The incident pressure and radial velocity are
\[
\begin{aligned}
	p^{\rm inc}(t,\rho)
	&=\frac{1}{2\pi}\int_0^\infty
	\widehat q(\kappa)J_0(\kappa\rho)
	G_c(\kappa,t)\,\kappa\,d\kappa,\\
	v_\rho^{\rm inc}(t,\rho)
	&=\frac{1}{2\pi}\int_0^\infty
	\widehat q(\kappa)J_1(\kappa\rho)
	G_s(\kappa,t)\,\kappa\,d\kappa.
\end{aligned}
\]
For $\rho>0$, the Cartesian incident velocity is
\[
v^{\rm inc}(t,x)
=v_\rho^{\rm inc}(t,\rho)\frac{x-x_s}{\rho}.
\]

The convolutions in $G_c$ and $G_s$ are evaluated analytically.

The primary incident table truncates the wavenumber integrals
at $K=60$. It has 4,096 quadrature nodes in $\kappa$ and
160 quadrature nodes for the source-radius integral defining
$\widehat q$. The table contains 1,801 radial nodes on
$0\leq\rho\leq18$ and 4,001 time nodes on $0\leq t\leq100$.
Its spacings are 
$
\Delta\rho=0.01,\qquad \Delta t_{\rm inc}=0.025.
$
Tensor-product cubic interpolation in $(t,\rho)$ supplies
the incident values at the requested points. We denote this
tabulated approximation by $\widetilde U^{\rm inc}_{60}$.

The scattered field satisfies the homogeneous acoustic equations
\[
\partial_t p^s+\nabla\cdot v^s=0,
\qquad
\partial_t v^s+\nabla p^s=0,
\qquad
U^s(0,x)=0.
\]
The rigid-wall condition is
\[
v^s\cdot\nu_j=-v^{\rm inc}\cdot\nu_j
\qquad\text{on }\partial D_j.
\]
In the PINN calculation, the right-hand side is supplied by
$\widetilde v^{\rm inc}_{60}$.
Both artificial boundary treatments are performed  on $U^s$.

\paragraph{\bf Computational and observation regions.}

Let
\[
B_R(0)=\{x\in\mathbb R^2:|x|<R\},
\qquad R=6.
\]
The common computational fluid region and observation region are
$$
\Omega_R=B_R(0)\cap\Omega,
\qquad
\Omega_{\rm obs}=\Omega_{\rm phys}=B_3(0)\cap\Omega.
$$
The filter buffer is the annulus

$$
\mathcal B=\{x:3<|x|<6\},
\qquad w=3,
$$
where $w$ denotes its width.
The Bayliss--Turkel condition is imposed on the same outer circle
$|x|=R$. On this circle, the outward unit normal of the
computational region is
$$
n_R(x)=\frac{x}{R}.
$$

\paragraph{\bf Second-order Bayliss--Turkel condition.}
We use polar coordinates about the origin,
\[
x=(r\cos\theta,r\sin\theta),
\qquad r=|x|,
\qquad 0\leq\theta<2\pi.
\]
  On $r=R$, it is easy to see that
\[
\partial_r
=\cos\theta\,\partial_{x_1}
+\sin\theta\,\partial_{x_2},
\qquad
\partial_\theta
=R\bigl(-\sin\theta\,\partial_{x_1}
+\cos\theta\,\partial_{x_2}\bigr).
\]

The second-order circular Bayliss--Turkel condition for the
scattered pressure is as follows
\begin{equation}\label{eq:test3-bt2}
	(\partial_t+R^{-1})
	\bigl(\partial_r+\partial_t+(2R)^{-1}\bigr)p^s
	=\frac{\partial_\theta^2p^s+p^s/4}{2R^2},
	\qquad r=R.
\end{equation}
The angular derivative supplies the second-order circular
correction in the radiation expansion of Bayliss and
Turkel (see~\cite{MR596431}).

The neural implementation introduces two scalar auxiliary fields,
$Q(t,\theta)$ and $W(t,\theta)$, defined on the outer circle by
\[
Q=\partial_rp^s+\partial_tp^s+\frac{p^s}{2R},
\qquad
W=\partial_\theta p^s.
\]
The resulting auxiliary system is
\begin{equation}\label{eq:test3-bt2-aux}
	\begin{aligned}
		\partial_rp^s+\partial_tp^s+\frac{p^s}{2R}-Q&=0,\\
		\partial_tQ+\frac{Q}{R}
		-\frac{\partial_\theta W+p^s/4}{2R^2}&=0,\\
		W-\partial_\theta p^s&=0.
	\end{aligned}
\end{equation}

Their initial values are
$$Q(0,\theta)=W(0,\theta)=0.$$

\paragraph{\bf Acoustic phase-space filter.}

Let $k=(k_1,k_2)$ denote the spatial Fourier wavevector.
The 2 propagating acoustic branches and their group velocities
are then
\[
\omega_\pm(k)=\pm|k|,
\qquad
v_{g,\pm}(k)=\nabla_k\omega_\pm(k)=\pm e(k),
\qquad
e(k)=\frac{k}{|k|},
\qquad k\neq0.
\]

For a Fourier coefficient
$\widehat V=(\widehat p,\widehat v_1,\widehat v_2)^T$,
we denote $\widehat v=(\widehat v_1,\widehat v_2)^T$.
The acoustic branch projectors are then
\[
P_\pm(k)\widehat V
=\frac12\bigl(\widehat p\pm e(k)\cdot\widehat v\bigr)
\begin{pmatrix}
	1\\
	\pm e_1(k)\\
	\pm e_2(k)
\end{pmatrix}.
\]

The filter-window index is $$\ell=0,\ldots,15.$$ 
We define the central angle and outward central normal of window
$\ell$ by
\[
\theta_\ell=\frac{2\pi\ell}{16},
\qquad
n_\ell=(\cos\theta_\ell,\sin\theta_\ell)^T,
\qquad
\Delta\theta=\frac{2\pi}{16}.
\]
Define the signed angular difference reduced modulo $2\pi$ by
\[
 d_\ell(\theta)
 =\bigl((\theta-\theta_\ell+\pi)\bmod 2\pi\bigr)-\pi
 \in[-\pi,\pi).
\]
Define the unnormalized angular factor
\[
\zeta_\ell(\theta)=
\begin{cases}
	\displaystyle
	\cos^2\!\left(\frac{\pi d_\ell(\theta)}
	{2\Delta\theta}\right),
	&|d_\ell(\theta)|<\Delta\theta,\\
	0,&|d_\ell(\theta)|\geq\Delta\theta,
\end{cases}
\]
and its normalized version
\[
\alpha_\ell(\theta)
=\frac{\zeta_\ell(\theta)}
{\left(\sum_{m=0}^{15}\zeta_m(\theta)^2\right)^{1/2}},
\]
which gives
\[
\sum_{\ell=0}^{15}\alpha_\ell(\theta)^2=1.
\]

The radial factor and full spatial window are
\[
\chi(r)=
\vartheta\!\left(\frac{r-3}{0.65w}\right)
\vartheta\!\left(\frac{R-r}{0.15w}\right),
\qquad
\eta_\ell(x)=\chi(|x|)\alpha_\ell(\theta).
\]

For $k\neq0$, the outgoing branch multipliers are
\[
m_{\ell,\pm}(k)
=\left[
1+\exp\!\left(
-\frac{\pm e(k)\cdot n_\ell-0.2}{0.05}
\right)
\right]^{-1}.
\]

Define the matrix-valued removal multiplier by
\[
M_\ell(k)
=m_{\ell,+}(k)P_+(k)+m_{\ell,-}(k)P_-(k),
\qquad k\neq0,
\qquad
M_\ell(0)=0.
\]

The filter uses a $192\times192$ Cartesian grid on
$[-R,R)^2$, with grid spacing
$$
h_{\rm filt}=\frac{2R}{192}=0.0625.
$$
Let $\mathcal F_h$ and $\mathcal F_h^{-1}$ denote the
two-dimensional discrete Fourier transform and its inverse,
applied separately to the three field components.
At a slab endpoint, let $V^{(0)}$ contain the raw neural
scattered field at grid points in the buffer annulus, with
zero values elsewhere on the Cartesian grid.
The windows are applied   in the order
$$\ell=0,\ldots,15$$ as follows

\[
V^{(\ell+1)}
=V^{(\ell)}
-\eta_\ell\operatorname{Re}\!\left\{
\mathcal F_h^{-1}
\left[M_\ell\,\mathcal F_h\bigl(\eta_\ell V^{(\ell)}\bigr)\right]
\right\}.
\]

We define the removed grid contribution by
\[
D_h=V^{(0)}-V^{(16)}.
\]
If $\mathcal I_hD_h$ denotes its componentwise cubic
interpolant, the filtered scattered endpoint at a fluid point
$x$ is
\[
U^{s,{\rm end}}(x)
=U^{s,{\rm raw}}(x)-(\mathcal I_hD_h)(x).
\]

Filtering is performed every
$\Delta t_{\rm filter}=0.25$.
The acoustic group speed is one, as such
$$
v_{\max}\Delta t_{\rm filter}
= 0.25<\frac{w}{3}=1.
$$
During each slab, the outer-circle closure for TDPSF is
\[
v^s\cdot n_R=0
\qquad\text{on }r=R.
\]

\paragraph{\bf Neural representation and input coordinates.}
Partition the simulation interval into 400 slabs as follows
\[
t_n=n\Delta t,\qquad
\Delta t=0.25,\qquad
n=0,\ldots,400.
\]
For $n=0,\ldots,399$, slab $n$ is $[t_n,t_{n+1}]$, with local time
\[
\tau=t-t_n,\qquad 0\leq\tau\leq\Delta t.
\]
The field network represents the three scattered components
$$(p^s,v_1^s,v_2^s).$$ 

For the network features, introduce the source center
$c_0=x_s$ in addition to the obstacle centers $c_1,c_2,c_3$.
 
Define the smoothed distances
$$
d_i(x)=\sqrt{|x-c_i|^2+\varepsilon^2},
\qquad \varepsilon=0.05.
$$

The field-network input vector consists of
\[
\begin{aligned}
	\Phi(t,x;t_n)=\Bigl(
	&\tau,\ x_1/R,\ x_2/R,\\
	&\{(d_i(x)-t)/R\}_{i=0}^3,\ q(x)/20,\\
	&\{\sin(\nu(d_i(x)-t))\}_{i=0,\ldots,3;\,\nu\in\{2,4.5\}},\\
	&\{\cos(\nu(d_i(x)-t))\}_{i=0,\ldots,3;\,\nu\in\{2,4.5\}}
	\Bigr).
\end{aligned}
\]

The field network has 2 hidden layers, each of width 96,
with activate function $$\sigma_f(z)=\sin(15z)$$ and three real outputs.
Including all weights and biases, its parameter count is
$12{,}003.
$
We denote its output as
$$
U_\lambda^s(t,x)
=\bigl(p_\lambda^s(t,x),
v_{1,\lambda}^s(t,x),
v_{2,\lambda}^s(t,x)\bigr)^T,
$$
where $\lambda$ collects its trainable parameters.

For BT2--PINN, we use a separate auxiliary network to represent
$(Q(t,\theta),W(t,\theta))$.
It has inputs
$
(\tau,\cos\theta,\sin\theta),
$
two hidden layers of width 24, activation
$$\sigma_a(z)=\sin(5z),$$  and 2 real outputs.
Its parameter count is
$746.$

\paragraph{\bf Collocation samples.}
In each slab, we use 2,048 PDE collocation points and 2,048
initial-matching points.

For each of these two sets, 683 spatial points are drawn
uniformly from $\Omega_{\rm obs}$ and 1,365 are drawn uniformly
from $\Omega_R$.

The PDE collocation times are sampled within the current slab, while
the initial-matching points are evaluated at $t=t_n$.

Each obstacle has 128 wall collocation samples, and the outer
circle has 192 boundary collocation samples.

Moreover, boundary positions are sampled uniformly in angle, with times
sampled within the slab.

In the case of BT2--PINN, we also use 192 outer-circle angles to match the
auxiliary fields at the beginning of the slab.

\paragraph{\bf Loss functions and their normalization.}
For any residual array $Z$, we define
\[
\operatorname{mse}(Z)
=\frac{1}{N_Z}\sum_{q=1}^{N_Z}|Z_q|^2,
\]
where $N_Z$ is the total number of scalar entries in the array,
including its sample and component indices.

In the formulas below, $p^s$ and $v^s$ denote the field-network
outputs, and $Q,W$ denote the  auxiliary-network outputs.
The three-component PDE residual is displayed as follows
\[
\mathcal R_{\rm PDE}
=
\begin{pmatrix}
	\partial_t p^s+\partial_{x_1}v_1^s+\partial_{x_2}v_2^s\\
	\partial_t v_1^s+\partial_{x_1}p^s\\
	\partial_t v_2^s+\partial_{x_2}p^s
\end{pmatrix}.
\]

Let $U_n^{s,{\rm end}}(x)$ be the scattered-field endpoint
transferred from the preceding slab.
At the beginning of slab $n$, define
\[
\delta U^s(x)
=U_\lambda^s(t_n,x)-U_n^{s,{\rm end}}(x).
\]
For $n=0$, the target is $U_0^{s,{\rm end}}=0$.
For later slabs, TDPSF transfers the filtered endpoint and
BT2 transfers the raw neural endpoint.

On obstacle $j$, the scalar rigid-wall residual is
\[
B_j(t,x)
=\bigl(v^s(t,x)+\widetilde v^{\rm inc}_{60}(t,x)\bigr)
\cdot\nu_j(x),
\qquad x\in\partial D_j.
\]
The complete objective is as thus
\begin{equation}\label{eq:test3-loss}
	\mathcal L
	=\operatorname{mse}(\mathcal R_{\rm PDE})
	+100\mathcal L_I
	+\frac53\sum_{j=1}^3\operatorname{mse}(B_j)
	+5\mathcal L_B
	+10^{-6}\mathcal R_{\rm out},
\end{equation}
where $\mathcal L_I$ is the initial-matching loss,
$\mathcal L_B$ is the outer-boundary loss, and
$\mathcal R_{\rm out}$ regularizes the final-layer weights
and biases.

For TDPSF, we set
\[
\mathcal L_I=\operatorname{mse}(\delta U^s),
\qquad
\mathcal L_B=\operatorname{mse}(v^s\cdot n_R)
\quad\text{on }r=R.
\]

In the case of BT2--PINN, letting $(Q_n(\theta),W_n(\theta))$ denote the auxiliary
endpoint transferred from the preceding slab, with
$Q_0=W_0=0$, we define the two-component initial mismatch
\[
\delta(Q,W)(\theta)
=
\begin{pmatrix}
	Q_\mu(t_n,\theta)-Q_n(\theta)\\
	W_\mu(t_n,\theta)-W_n(\theta)
\end{pmatrix},
\]
as well as the initial loss  
\[
\mathcal L_I
=\operatorname{mse}(\delta U^s)
+\frac14\operatorname{mse}(\delta(Q,W)),
\]
the boundary residuals 
\[
\begin{aligned}
	b_1&=\partial_rp^s+\partial_tp^s+\frac{p^s}{2R}-Q,\\
	b_2&=\partial_tQ+\frac QR
	-\frac{\partial_\theta W+p^s/4}{2R^2},\\
	b_3&=\frac{W-\partial_\theta p^s}{R},
\end{aligned}
\]
and the outer-boundary loss  
\[
\mathcal L_B
=\frac{
	\operatorname{mse}(b_1)
	+\operatorname{mse}(b_2)
	+\operatorname{mse}(b_3)}{3}.
\]

Let $\mathsf A_f$ and $\mathsf b_f$ denote the weight matrix
and bias vector of the field network's final linear layer.
For BT2--PINN, let $\mathsf A_a$ and $\mathsf b_a$ denote the
corresponding auxiliary-network quantities.
Then, we have
\[
\mathcal R_{\rm out}
=
\begin{cases}
	\|\mathsf A_f\|_F^2+\|\mathsf b_f\|_2^2,
	&\text{ for TDPSF},\\
	\|\mathsf A_f\|_F^2+\|\mathsf b_f\|_2^2
	+\|\mathsf A_a\|_F^2+\|\mathsf b_a\|_2^2,
	&\text{ with BT2}.
\end{cases}
\]
The squared Frobenius norm $\|\mathsf A\|_F^2$ is the sum
of the squared matrix entries, and $\|\mathsf b\|_2^2$ is
the sum of the squared vector entries.

\paragraph{\bf Optimization and endpoint transfer.}
Each slab begins with a linear least-squares update of the
output weights and biases, while the hidden-layer parameters are
fixed. 

This is followed by all-weight limited-memory
Broyden--Fletcher--Goldfarb--Shanno (L-BFGS) iterations and then
another output least-squares update. 

The first slab has at most
1,800 L-BFGS iterations.

Subsequent slabs with left endpoint
$t_n<16$ have at most 600 iterations, and slabs with $t_n\geq16$
have at most 150 iterations. Both boundary methods use this same
training schedule.

In the linear solves, we use relative rank tolerance $10^{-5}$.

Initial matching transfers the preceding computed field
endpoint, including the filter correction to the filter TDPSF.
Algorithm BT2 also transfers its computed $Q$ and $W$ endpoints.
 
The 3 paired seeds, 0, 1, and 2, have the same field architecture,
collocation rules, and iteration caps for the two boundary
methods.

At $t=100$, we evaluate the saved field and auxiliary weights with
local slab origin $99.75$.

For TDPSF--PINN, we reconstruct the
annular-grid removal and subtract its interpolation from the neural
scattered field. We then add the prescribed incident field to form
the total field. 

\paragraph{\bf Training refinement during the forcing interval.}

We assess training resolution on $0\leq t\leq8$ with the same
three paired seeds $0,1,2$ for both methods.

In the joint refinement, the two field
layers have width $128$, the BT2 auxiliary layers have width $32$,
and each slab uses $4096$ PDE points and $4096$ initial points.

Each obstacle has $256$ wall samples, and the outer boundary and
BT2 auxiliary initial trace each have $384$ samples.

The first slab
has at most $3600$ L-BFGS iterations and each later slab at most
$1200$.

The field and BT2 auxiliary parameter counts are $20099$
and $1250$, respectively.

The physical model, source, geometry,
reference, slab width, filtering grid, loss weights, ridge coefficient,
and linear-solve rank tolerance are the same in both configurations.

This joint comparison measures the combined effect of increased
network capacity, sampling, and optimization budgets.

At endpoints $$t=4,5,6,8$$ we save the trained networks and independently
sample the preceding quarter-unit slab.

We use NumPy's independent
generator with seed $731207$ and the slab index. 

Each PDE sample set
contains $2730$ spatial points in the observation fluid region and
$5462$ in the full computational fluid region, with independent
uniform times on that slab.

The initial-trace sample set has the
same spatial counts at the slab's left endpoint.

We additionally
use $1024$ independent space--time samples on each obstacle and
$2048$ on the outer circle. 

BT2 auxiliary initial matching has
$2048$ independent angles. These validation coordinates are common
to both methods, both configurations, and all 3 seeds.

We evaluate the component-mean squared PDE residual, the field
initial mismatch, each obstacle's normal-velocity residual, and
each outer-boundary residual. In the case of BT2--PINN, we record $b_1,b_2,b_3$
separately and also evaluate the initial mismatch of $Q$ and $W$.

\paragraph{\bf Observation quadrature, error, and energy.}
The observation quadrature starts from a $97\times97$
Cartesian grid on $[-3,3]^2$ and retains the points in
$\Omega_{\rm obs}$. Denote these retained points by
\[
x_q^{\rm obs},\qquad q=1,\ldots,N_{\rm obs},
\qquad N_{\rm obs}=5{,}875.
\]
The grid spacing and cell-area weight are
$$
h_{\rm obs}= 0.0625,
 w_{\rm obs}=h_{\rm obs}^2.$$
For a 3-component field $F$, we define the discrete
observation norm
\[
\|F\|_{{\rm obs},h}
=\left(
w_{\rm obs}\sum_{q=1}^{N_{\rm obs}}
|F(x_q^{\rm obs})|^2
\right)^{1/2},
\]
where $|F|^2$ is the sum of the squared components.

Let $U_h^s(t_n)$ denote the computed scattered endpoint:
the filtered endpoint for TDPSF and the raw neural endpoint
for BT2.  
Let $U_{\rm ref}^s(t_n)$ denote the independent exterior
scattered-field reference.
The total fields used in the comparison are as thus
\[
U_h(t_n)
=\widetilde U^{\rm inc}_{60}(t_n)+U_h^s(t_n),
\qquad
U_{\rm ref}(t_n)
=\widetilde U^{\rm inc}_{60}(t_n)+U_{\rm ref}^s(t_n).
\]
Define the single reference scale
\[
S_{\rm ref}
=\max_{0\leq m\leq400}
\|U_{\rm ref}(t_m)\|_{{\rm obs},h}.
\]

We set  $U_h=(p_h,v_{1,h},v_{2,h})^T$ and
$v_h=(v_{1,h},v_{2,h})^T$.
The physical-region and computational-region energies are
\[
\mathcal E_{\rm phys}(t)
=\frac12\int_{\Omega_{\rm obs}}
\bigl(p_h^2+|v_h|^2\bigr)\,dx,
\qquad
\mathcal E_{\rm comp}(t)
=\frac12\int_{\Omega_R}
\bigl(p_h^2+|v_h|^2\bigr)\,dx.
\]
These integrals are evaluated using the corresponding
Cartesian cell-area quadratures.

The trajectory data also record total pressure at the
4 receiver positions
\[
x_1^{\rm rec}=(-2,1),\qquad
x_2^{\rm rec}=(2,1),\qquad
x_3^{\rm rec}=(0,2.7),\qquad
x_4^{\rm rec}=(0,-2.5).
\] 

\paragraph{\bf Independent exterior reference.}
The reference represents the scattered pressure using
outgoing cylindrical harmonics centered at the three
obstacles. We follow the circular multipole approach for
rigid-cylinder scattering~\cite{Sherer2004}.
For obstacle $j$, let $\phi_j$ denote the polar angle about
$c_j$. Angular order $M$ retains the constant angular mode
and the functions
\[
\cos(m\phi_j),\qquad \sin(m\phi_j),
\qquad m=1,\ldots,M.
\]
The baseline reference have $M=28$.

Define
\[
\widehat p^s(s,x)
=\int_0^\infty e^{-st}p^s(t,x)\,dt.
\]
 
The scattered velocity follows as 
\[
\widehat v^s(s,x)
=-\frac{1}{s}\nabla\widehat p^s(s,x).
\]

The numerical inversion samples
$$
s=\gamma_{\rm L}-\mathrm i\omega,
$$
where $\mathrm i^2=-1$, $\gamma_{\rm L}>0$ is the Laplace damping
parameter, and $\omega$ is temporal angular frequency.
We use the parameters
\[
\gamma_{\rm L}=0.05,\qquad
\omega_{\rm cut}=32,\qquad
T_{\rm inv}=256,\qquad
N_{\rm inv}=4{,}096,
\]
where $\omega_{\rm cut}$ is the retained temporal-frequency
cutoff, $T_{\rm inv}$ is the period parameter of the
discrete inversion, and $N_{\rm inv}$ is its number of
temporal nodes. They give
\[
\Delta\omega=\frac{2\pi}{T_{\rm inv}},
\qquad
\Delta t_{\rm inv}
=\frac{T_{\rm inv}}{N_{\rm inv}}=0.0625.
\]

Reference sensitivity is evaluated separately with angular order
$M=32$ and damping $\gamma_{\rm L}=0.07$, keeping the other
parameters fixed.
For a varied reference $U_{\rm ref}^{\rm var}$, define
\[
D_{\rm ref}(t_n)
=\frac{
\|U_{\rm ref}^{\rm var}(t_n)-U_{\rm ref}(t_n)\|_{{\rm obs},h}}
{S_{\rm ref}}.
\]
The recorded differences at selected times are displayed in the table below
\[
\begin{array}{lccc}
\text{Changed parameter}&D_{\rm ref}(2)&D_{\rm ref}(4)&D_{\rm ref}(100)\\
\hline
M=32&3.34\times10^{-7}&1.59\times10^{-6}&3.08\times10^{-11}\\
\gamma_{\rm L}=0.07&9.55\times10^{-11}&9.28\times10^{-11}&6.02\times10^{-9}
\end{array}
\]

\paragraph{\bf Incident-field resolution.}
The total-field differences satisfy
\[
U_h-U_{\rm ref}=U_h^s-U_{\rm ref}^s,
\qquad
U_{\rm ref}^{\rm var}-U_{\rm ref}
=U_{\rm ref}^{s,{\rm var}}-U_{\rm ref}^s.
\]
The PINN wall forcing is evaluated from
$\widetilde U^{\rm inc}_{60}$, while the independent scattered
reference uses the Laplace-transformed source directly.

Let $U_K^{\rm inc}$ denote direct numerical evaluation of the
incident Fourier--Bessel integrals with radial-wavenumber cutoff
$K$. The cutoff comparisons have $$K\in\{60,120,240\},$$ 
16,384 wavenumber quadrature nodes, and 256 source-radius
quadrature nodes.

We evaluate the fields at the 5,875 observation
points and on a wall set $\mathcal W$ containing 128 points on
each obstacle, for 384 wall points in total. The wall points on
obstacle $j$ are
\[
c_j+a(\cos\varphi_m,\sin\varphi_m),
\qquad
\varphi_m=\frac{2\pi(m+0.314159)}{128},
\qquad m=0,\ldots,127.
\]
For $$\xi\in\mathcal W\cap\partial D_j,$$ set
$$\nu(\xi)=\nu_j(\xi).$$  The 9 evaluation times are
\[
\mathcal T_{\rm inc}
=\{1.13,2,3.73,4.5,6.25,8,12,20,100\}.
\]
For 2 cutoffs $K_a$ and $K_b$, define
\[
D_{\rm inc}(K_a,K_b;t)
=\frac{\|U_{K_a}^{\rm inc}(t)-U_{K_b}^{\rm inc}(t)\|_{{\rm obs},h}}
{S_{\rm ref}},
\]
and the absolute wall-forcing difference
\[
D_{\rm wall}(K_a,K_b;t)
=\max_{\xi\in\mathcal W}
\left|
\bigl(v_{K_a}^{\rm inc}(t,\xi)-v_{K_b}^{\rm inc}(t,\xi)\bigr)
\cdot\nu(\xi)
\right|.
\]
Representative values during the forcing interval and at the
final time are given below
\[
\begin{array}{cccc}
t&(K_a,K_b)&D_{\rm inc}(K_a,K_b;t)&D_{\rm wall}(K_a,K_b;t)\\
\hline
2&(60,120)&3.38\times10^{-4}&1.42\times10^{-4}\\
2&(120,240)&5.44\times10^{-5}&2.49\times10^{-5}\\
4.5&(60,120)&2.05\times10^{-3}&8.73\times10^{-4}\\
4.5&(120,240)&3.31\times10^{-4}&1.52\times10^{-4}\\
100&(60,120)&6.13\times10^{-13}&1.32\times10^{-13}\\
100&(120,240)&1.32\times10^{-12}&2.03\times10^{-13}
\end{array}
\]

At fixed cutoff $K=60$, we also compare the tabulated incident
field, with spacings $\Delta\rho=0.01$ and
$\Delta t_{\rm inc}=0.025$, with the direct evaluation using
16,384 wavenumber nodes and 256 source-radius nodes. The table below
reports the normalized observation-norm difference
\[
\frac{\|\widetilde U^{\rm inc}_{60}(t)-U^{\rm inc}_{60}(t)\|_{{\rm obs},h}}
{S_{\rm ref}}
\]
and the absolute maximum difference of the normal wall velocity
on the same set $\mathcal W$:
\[
\begin{array}{ccc}
t&\text{Incident difference}&\text{Wall-normal difference}\\
\hline
2&1.09\times10^{-6}&4.18\times10^{-7}\\
4.5&6.56\times10^{-6}&2.55\times10^{-6}\\
100&2.21\times10^{-13}&4.94\times10^{-14}
\end{array}
\]

An independent 2,048-node Gauss--Legendre quadrature in time
checks the analytic convolution. 

With $K=60$, 4,096 wavenumber
nodes, and 160 source-radius nodes, the resulting normalized
incident-field differences are $1.38\times10^{-13}$ at $t=2$
and $1.99\times10^{-12}$ at $t=4.5$, using the same
observation norm and scale $S_{\rm ref}$.

\subsubsection{Numerical results for Test 3}
Figure~\ref{fig:test3-accuracy} displays the complete fixed-scale
error histories through $T=100$, including the forcing interval
$0<t<8$.

At $t=5$, the median
errors are $0.149$ for TDPSF--PINN and $0.160$ for BT2--PINN;
at $t=6$, they are $0.162$ and $0.169$, respectively.

The late-time histories demonstrate the sustained accuracy of
TDPSF--PINN: its median error is $5.97\times10^{-3}$ at $t=32$
and $8.68\times10^{-8}$ at $t=100$, compared with
$6.90\times10^{-2}$ and $9.79\times10^{-2}$ for BT2--PINN.

For the slab ending at $t=8$, the median independent PDE mean squares are $2.12\times10^{-4}$ and $8.33\times10^{-5}$ for the baseline and refined TDPSF--PINN calculations, and $2.22\times10^{-4}$ and $8.67\times10^{-5}$ for BT2--PINN, respectively.

Figure~\ref{fig:test3-energy} demonstrates the corresponding suppression
of residual energy by TDPSF--PINN in the physical region, buffer
annulus, and entire computational domain.

At $T=100$, the median
physical-region and buffer energies are $1.24\times10^{-11}$
and $4.63\times10^{-11}$. The median total computational energies
are $5.87\times10^{-11}$ for TDPSF--PINN and
$2.88\times10^{-2}$ for BT2--PINN.

Together with the error
histories, these results show the long-time advantage of
TDPSF--PINN in this  scattering problem.

All 6 primary
trajectories reach $T=100$ with finite computed fields.

For seeds 0, 1, and 2, the complete trajectories use
$2606.103$, $2430.155$, and $2488.520$ process CPU seconds
for TDPSF--PINN and $4567.839$, $4573.288$, and $4565.992$
for BT2--PINN, respectively. These costs include training and
endpoint operations under the shared iteration schedule.

\begin{figure}[htbp]\centering\color{black}\captionsetup{font={color=black}}
\includegraphics[width=\linewidth]{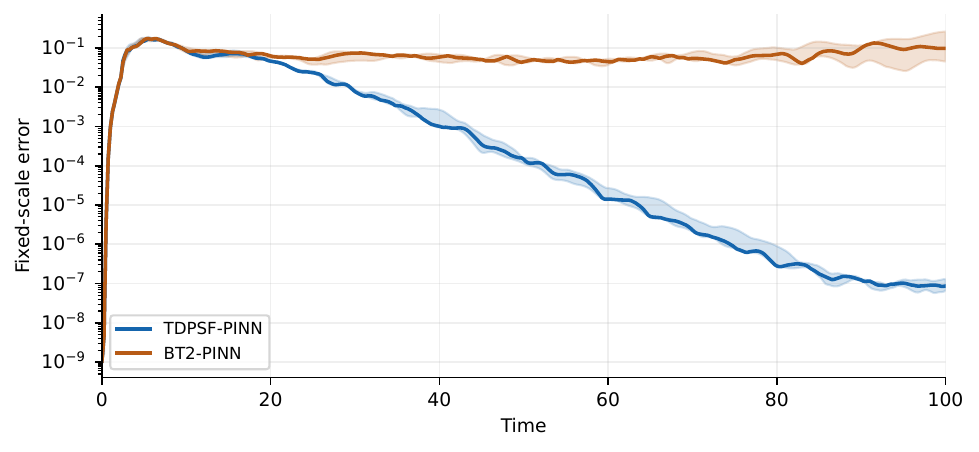}
\caption{Test 3: fixed-scale error histories for TDPSF--PINN and
BT2--PINN through $T=100$. Curves are pointwise medians and bands
show the full range over three paired seeds. TDPSF--PINN achieves
sustained late-time error decay.}
\label{fig:test3-accuracy}
\end{figure}

\begin{figure}[htbp]\centering\color{black}\captionsetup{font={color=black}}
\includegraphics[width=\linewidth]{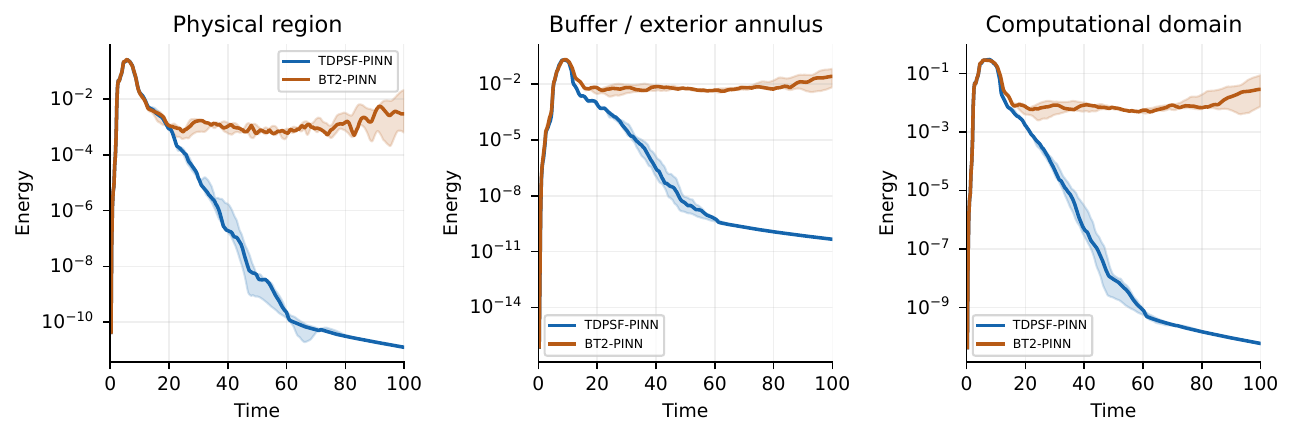}
\caption{Test 3: physical-region, buffer-annulus, and computational-domain
energy histories through $T=100$. Curves are pointwise medians and
bands show the full range of the PINN total-field energies over
the same three paired seeds.}
\label{fig:test3-energy}
\end{figure}

\space
\FloatBarrier
\end{revision}

\Needspace{8\baselineskip}
\section{Conclusion}
\label{sec:conclusion}
\begin{revision}
In our work, we present a microlocal PINN Filter (TDPSF--PINN), that  couples a slabwise neural residual solver with the
established TDPSF  mechanism. 
The 3   tests evaluate this coupling for $3$ problems: anisotropic focusing
NLS, zero-flow acoustic radiation, and scattering by three rigid disks. In the $3$ tests,   TDPSF--PINN is compared with  3 methods: PML, CRBC and BT2 (all coupled with PINN).
The comparisons are done for long time windows.
In the 3 reported
long-time simulations, the  TDPSF--PINN method is more stable and achieves smaller  errors in all three
tests.

\section{Statement of AI use}\label{sec:ai}

	ChatGPT (OpenAI) was used to assist with drafting, debugging the Python implementation and preparing numerical illustrations. The authors reviewed the code and take full responsibility for the content.   	The authors also used 	ChatGPT   for  manuscript polishing and checking.  The authors have checked and ensured that everything is correct.
	 
\end{revision}
\bibliographystyle{amsplain}

\begin{thebibliography}{99}
\setlength{\baselineskip}{12pt}
\setlength{\lineskiplimit}{0pt}
\setlength{\lineskip}{1pt}

\bibitem{abarbanel:PMLinstability2}
S. Abarbanel and D. Gottlieb,
\newblock On the construction and analysis of absorbing layers in CEM.
\newblock {\em Applied Numerical Mathematics}, 27(4):331--340, 1998.

\bibitem{abarbanel:PMLinstability}
S. Abarbanel, D. Gottlieb, and J. S. Hesthaven,
\newblock Long Time Behavior of the Perfectly Matched Layer Equations in Computational Electromagnetics.
\newblock {\em Journal of Scientific Computing}, 17(1):405--422, 2002.

\bibitem{AbediPardoAlkhalifah2025}
M.~M. Abedi, D.~Pardo, and T.~Alkhalifah,
\newblock Gabor-enhanced physics-informed neural networks for fast simulations of acoustic wavefields.
\newblock {\em Neural Networks}, 193 (2026). DOI: 10.1016/J.NEUNET.2025.107978.

\bibitem{AlkhadhrAlmekkawy2023}
S.~Alkhadhr and M.~Almekkawy,
\newblock Wave Equation Modeling via Physics-Informed Neural Networks: Models of Soft and Hard Constraints for Initial and Boundary Conditions.
\newblock {\em Sensors}, 23(5):2792, 2023.

\bibitem{AlpertGreengardHagstrom2000}
\begingroup\color{black}
B. Alpert, L. Greengard, and T. Hagstrom,
\newblock Rapid Evaluation of Nonreflecting Boundary Kernels for Time-Domain Wave Propagation.
\newblock {\em SIAM Journal on Numerical Analysis}, 37(4):1138--1164, 2000.
\par\endgroup

\bibitem{MR596431}
A. Bayliss and E. Turkel,
\newblock Radiation Boundary Condition for Wave-Like Equations.
\newblock {\em Communications on Pure and Applied Mathematics}, 33(6):707--725, 1980.

\bibitem{BecacheFauqueuxJoly2003}
E. B\'ecache, S. Fauqueux, and P. Joly,
\newblock Stability of perfectly matched layers, group velocities and anisotropic waves.
\newblock {\em Journal of Computational Physics}, 188(2):399--433, 2003.

\bibitem{Berenger1994}
J.-P. B\'erenger,
\newblock A Perfectly Matched Layer for the Absorption of Electromagnetic Waves.
\newblock {\em Journal of Computational Physics}, 114:185--200, 1994.

\bibitem{MR1412240}
J.-P. B{\'e}renger,
\newblock Three-dimensional Perfectly Matched Layer for the Absorption of Electromagnetic Waves.
\newblock {\em Journal of Computational Physics}, 127(2):363--379, 1996.

\bibitem{BerroneCanutoPintoreSukumar2023}
S.~Berrone, C.~Canuto, M.~Pintore, and N.~Sukumar,
\newblock Enforcing Dirichlet boundary conditions in physics-informed neural networks and variational physics-informed neural networks.
\newblock {\em Heliyon}, 9(8):e18820, 2023.

\bibitem{Daubechies1988}
I. Daubechies,
\newblock Time-frequency localization operators: A geometric phase space approach.
\newblock {\em IEEE Transactions on Information Theory}, 34(4):605--612, 1988.

\bibitem{Dohnal2009}
\begingroup\color{black}\hypersetup{urlcolor=red}
T. Dohnal,
\newblock Perfectly matched layers for coupled nonlinear Schr\"odinger equations with mixed derivatives.
\newblock {\em Journal of Computational Physics}, 228(23):8752--8765, 2009.\newline DOI: \url{https://doi.org/10.1016/j.jcp.2009.08.023}.
\newblock \url{https://arxiv.org/abs/0905.2321}.
\par\endgroup

\bibitem{MR0471386}
B. Engquist and A. Majda,
\newblock Absorbing Boundary Conditions for the Numerical Simulation of Waves.
\newblock {\em Proceedings of the National Academy of Sciences of the United States of America}, 74(5):1765--1766, 1977.

\bibitem{MR0436612}
B. Engquist and A. Majda,
\newblock Absorbing Boundary Conditions for the Numerical Simulation of Waves.
\newblock {\em Mathematics of Computation}, 31(139):629--651, 1977.

\bibitem{MR517938}
B. Engquist and A. Majda,
\newblock Radiation boundary conditions for acoustic and elastic wave calculations.
\newblock {\em Communications on Pure and Applied Mathematics}, 32(3):313--357, 1979.



\bibitem{givoli1991non}
\begingroup\color{black}\hypersetup{urlcolor=red}D. Givoli,
\newblock Non-reflecting boundary conditions.
\newblock {\em Journal of computational physics}, 94.1 (1991): 1-29.\endgroup

\bibitem{MR1819643}
T. Hagstrom,
\newblock Radiation Boundary Conditions for the Numerical Simulation of Waves.
\newblock {\em Acta Numerica}, 8:47--106, 1999.

\bibitem{MR2032866}
T. Hagstrom,
\newblock New Results on Absorbing Layers and Radiation Boundary Conditions.
\newblock {\em Topics in Computational Wave Propagation}, volume~31 of {\em Lecture Notes in Computational Science and Engineering}, pages 1--42. Springer, Berlin, 2003.

\bibitem{HagstromNguyenSofferStucchioTran2026}
T. Hagstrom, D. P. C. Nguyen, A. Soffer, C. Stucchio, and M.-B. Tran,
\newblock A Time-Dependent Phase Space Filter for Anisotropic Wave Equations on Unbounded Domains.
\newblock {\em SIAM Journal on Scientific Computing}, 48(2):A804--A827, 2026. DOI: 10.1137/25M1748536.

\bibitem{HagstromWarburton2009}
\begingroup\color{black}\hypersetup{urlcolor=red}
T. Hagstrom and T. Warburton,
\newblock Complete radiation boundary conditions: minimizing the long time error growth of local methods.
\newblock {\em SIAM Journal on Numerical Analysis}, 47(5):3678--3704, 2009.
\newblock \url{https://doi.org/10.1137/090745477}.
\par\endgroup

\bibitem{hu:unstablePML}
F. Q. Hu,
\newblock On Absorbing Boundary Conditions for Linearized Euler Equations by a Perfectly Matched Layer.
\newblock {\em Journal of Computational Physics}, 129(1):201--219, 1996.

\bibitem{JiangGreengard2008}
S. Jiang and L. Greengard,
\newblock Efficient representation of nonreflecting boundary conditions for the time-dependent Schr\"odinger equation in two dimensions.
\newblock {\em Communications on Pure and Applied Mathematics}, 61(2):261--288, 2008.

\bibitem{RaissiPerdikarisKarniadakis2019}
M. Raissi, P. Perdikaris, and G. E. Karniadakis,
\newblock Physics-Informed Neural Networks: A Deep Learning Framework for Solving Forward and Inverse Problems Involving Nonlinear Partial Differential Equations.
\newblock {\em Journal of Computational Physics}, 378:686--707, 2019.

\bibitem{Rauch2012}
J. Rauch,
\newblock {\em Hyperbolic Partial Differential Equations and Geometric Optics}, Graduate Studies in Mathematics, Vol. 133, American Mathematical Society, Providence, RI, 2012.

\bibitem{RenRaoChenWangSunLiu2024}
\begingroup\color{black}
P.~Ren, C.~Rao, S.~Chen, J.-X. Wang, H.~Sun, and Y.~Liu,
\newblock SeismicNet: Physics-informed neural networks for seismic wave modeling in semi-infinite domain.
\newblock {\em Computer Physics Communications}, 295:109010, 2024.
\par\endgroup

\bibitem{SofferStucchio2007}
A. Soffer and C. Stucchio,
\newblock Open boundaries for the nonlinear Schr\"odinger equation.
\newblock {\em Journal of Computational Physics}, 225(2):1218--1232, 2007. DOI: 10.1016/j.jcp.2007.01.020.

\bibitem{us:multiscale}
A. Soffer and C. Stucchio,
\newblock Multiscale Resolution of Shortwave-Longwave Interactions.
\newblock {\em Communications on Pure and Applied Mathematics}, 62(1):82--124, 2009.

\bibitem{SofferStucchioTran2023}
A. Soffer, C. Stucchio, and M.-B. Tran,
\newblock {Time Dependent Phase Space Filters: A Stable Absorbing Boundary Condition}, Springer Briefs in PDEs and Data Science, Springer Nature, Singapore, 2023.

\bibitem{SukumarSrivastava2022}
\begingroup\color{black}
N.~Sukumar and A.~Srivastava,
\newblock Exact imposition of boundary conditions with distance functions in physics-informed deep neural networks.
\newblock {\em Computer Methods in Applied Mechanics and Engineering}, 389:114333, 2022.
\par\endgroup

\bibitem{WuAghamiryOpertoMa2022}
Y.~Wu, H.~S. Aghamiry, S.~Operto, and J.~Ma,
\newblock Wave simulation in non-smooth media by PINN with quadratic neural network and PML condition.
\newblock {\em arXiv preprint arXiv:2208.08276}, 2022.

\bibitem{Zheng2007}
\begingroup\color{black}\hypersetup{urlcolor=black}
C. Zheng,
\newblock A perfectly matched layer approach to the nonlinear Schr\"odinger wave equations.
\newblock {\em Journal of Computational Physics}, 227(1):537--556, 2007.
\newblock \url{https://doi.org/10.1016/j.jcp.2007.08.004}.
\par\endgroup

\bibitem{Sherer2004}
\begingroup\color{black}\hypersetup{urlcolor=black}
S. E. Sherer,
Scattering of sound from axisymmetric sources by multiple circular cylinders,
{\em Journal of the Acoustical Society of America}, 115(2):488--496, 2004.
DOI: \url{https://doi.org/10.1121/1.1641790}.
\par\endgroup

\end{thebibliography}

\end{document}